\documentclass[11pt]{article}
\usepackage[margin=1in]{geometry}
\usepackage[T1]{fontenc}
\usepackage{lmodern}
\usepackage{amsmath,amssymb,amsthm,mathtools,bm,booktabs,array,microtype}
\usepackage[unicode,hidelinks,bookmarks=false]{hyperref}
\hypersetup{
  pdftitle={Optimal control of quasi-static poroviscoelasticity with vanishing storage and structural viscosity},
  pdfauthor={Jeonghun J. Lee},
  pdfkeywords={Biot consolidation, Kelvin--Voigt viscoelasticity, optimal control, singular limits}
}
\usepackage{enumitem}
\numberwithin{equation}{section}
\newtheorem{theorem}{Theorem}[section]
\newtheorem{lemma}[theorem]{Lemma}
\newtheorem{proposition}[theorem]{Proposition}
\newtheorem{corollary}[theorem]{Corollary}
\newtheorem{assumption}[theorem]{Assumption}
\theoremstyle{definition}
\theoremstyle{remark}\newtheorem{remark}[theorem]{Remark}
\newcommand{\eps}{\varepsilon}
\newcommand{\norm}[1]{\lVert#1\rVert}
\newcommand{\abs}[1]{\lvert#1\rvert}

\newcommand{\dual}[2]{\langle #1,#2\rangle}
\newcommand{\dd}{\,\mathrm d}
\newcommand{\UU}{\mathcal U}
\newcommand{\YY}{\mathcal Y}
\newcommand{\KK}{\mathcal K}
\newcommand{\ee}{\mathrm e}
\newcommand{\vv}{\mathrm v}
\newcommand{\ad}{\mathrm{ad}}
\newcommand{\calL}{\mathcal L}
\DeclareMathOperator{\divergence}{div}
\title{Optimal control of quasi-static poroviscoelasticity with vanishing storage and structural viscosity}
\author{Jeonghun J. Lee\\
Department of Mathematics, Baylor University \\One Bear Place \#97328, Waco, TX 76798, USA\\
\texttt{jeonghun\_lee@baylor.edu}}
\date{}
\begin{document}
\maketitle

\begin{abstract}

We study linear-quadratic optimal control problems governed by a two-field quasi-static poroviscoelastic system with constrained storage $c_0\geq0$ and Kelvin--Voigt coefficient $\delta\geq0$. 
We formulate the state and adjoint equations in the four parameter regimes determined by positivity of $c_0$ and $\delta$, identify the corresponding initial and terminal conditions, and derive the optimality conditions. Under a mechanical boundary condition that ensures a pressure-to-dilation inf--sup condition, for mechanical loads and fluid sources in suitable function spaces and compatible initial data, we prove that, as both coefficients vanish, the states, adjoint, control variables converge to the solution of the quasi-static Biot model without storage at certain rates depending on the storage and the viscosity coefficients. We also present some explicit examples which show that the convergence rates in our analyses are sharp. 
\end{abstract}



\section{Introduction}
The quasi-static Biot model \cite{Biot1941} describes the deformation of a fluid-saturated porous medium. Let $u$ denote the displacement of the porous skeleton and let $p$ denote the fluid pressure. With a Kelvin--Voigt constitutive law and the classical fluid-content relation, the governing equations are
\begin{align}
-\divergence\!\left(\mathcal C_{\rm e}\eps(u)+\delta\mathcal C_{\rm v}\eps(u_t)-\alpha pI\right)&=f,
\label{intro:momentum}\\
\partial_t\bigl(c_0p+\alpha\divergence u\bigr)-\divergence(\mathbf K\nabla p)&=g.
\label{intro:fluid}
\end{align}
Here $\eps(u)$ is the symmetric gradient of $u$, $\mathcal C_{\rm e}$ and $\mathcal C_{\rm v}$ are the elasticity and viscosity tensors, $\alpha>0$ is the Biot--Willis coefficient, $\mathbf K$ is the permeability tensor, $f$ is a body force, $g$ is a fluid source, $c_0\geq0$ is the constrained storage coefficient, and $\delta\geq0$ is the structural viscosity coefficient. When either coefficient vanishes, the character of the system changes: the admissible data, the independent initial conditions, and the terminal conditions of the adjoint problem depend on the parameter regime. We study linear-quadratic optimal control problems for this system and their behavior in the joint limit $c_0,\delta\to0$, which leads to the quasi-static Biot model without storage.

When storage and structural viscosity are small, a natural question is whether an optimal control computed for the simpler poroelastic model without storage is a good approximation of the optimal control for the full poroviscoelastic model, and at which rate. We show that the answer depends not only on the material parameters but also on where the control acts: distributed fluid sources and tractions lead to first-order convergence in the viscosity, whereas normal Darcy fluxes lead to the rate $\delta^{3/4}$, and this difference is intrinsic for the corresponding control-to-state operators.

Existence and uniqueness for Biot's consolidation model were proved in \cite{Zenisek1984} by a compactness method. The forward problem also fits into the theory of degenerate evolution equations \cite{Showalter1974,Showalter2000}, and its well-posedness in the different storage--viscosity regimes was analyzed in \cite{BGSSW2016,BMW2023}, and an inf--sup theory for the quasi-static Biot equations with stability estimates robust in the material parameters was developed in \cite{KZ2026}; see also \cite{Verri2018,BGSV2019} for the roles of structural viscosity and compressibility in biomechanical models. Optimal control of the poroelastic problem without storage was studied in \cite{BS2022}, and of the poroviscoelastic problem in \cite{BSvisco2022,Strikwerda2023}; integer-valued controls with total-variation regularization for both systems were considered in \cite{BMSS2024}. Numerical methods for related control problems were developed in \cite{KLS2026,GMV2026,LeeThermo2026,LeeLocking2026}. Our main interest is in quantitative estimates for the dependence of the state, adjoint, and optimal control on the storage and structural viscosity coefficients, including their simultaneous vanishing; the time-dependent permeability considered in \cite{BSvisco2022,Strikwerda2023} is admitted in these estimates. Under our standing assumption that the skeleton is not fully clamped, the derivation of the optimality system requires no additional time-regularity assumptions on the pressure tracking data or the optimal pressure, in contrast to \cite{BS2022}; and our estimates do not require the elastic and viscous bilinear forms to be proportional.

Our main results are as follows.
\begin{enumerate}[label=(\roman*),leftmargin=*]
	\item For mechanical loads that are $H^1$ in time, $L^2$ fluid sources, and compatible initial data, we prove that, as $c_0,\delta\to0$, the states converge to the solution of the degenerate problem with $c_0=\delta=0$ at the rate $O(c_0+\sqrt\delta)$ in the energy norm and at the rate $O(c_0+\delta)$ in a weaker tracking norm by a duality argument. Both exponents of $\delta$ are sharp; for the pressure this holds already for fluid-source data.
	\item We prove that, as $c_0,\delta\to0$, the optimal controls converge to the optimal control of the degenerate problem with $c_0=\delta=0$, at the rate $O(c_0+\delta)$ for distributed fluid-source controls and for traction controls that are $H^1$ in time with zero initial trace, and at the rate $O(c_0+\delta^{3/4})$ for normal Darcy flux controls. An explicit two-dimensional example shows that the exponent $3/4$ cannot be improved in general, both for the flux-to-state operators and for the map from the target to the unconstrained optimal flux control.
	\item For fixed observation residuals, the adjoint variables are bounded uniformly in $c_0$ and $\delta$. The mechanical adjoint converges in $L^2(0,T;H^1(\Omega)^d)$, and the pressure adjoint converges at the rate $O(c_0+\delta)$ in $L^2(0,T;L^2(\Omega))$ and $O(c_0+\sqrt\delta)$ in $L^2(0,T;H^1(\Omega))$. In $L^2(0,T;L^2(\Omega))$, the optimal pressure adjoints converge at the rate of the corresponding optimal states.
	\item These results rest on a formulation of the state and adjoint equations in all four parameter regimes without additional time-differentiability assumptions on the tracking data. We identify the independent initial data and the terminal conditions in each regime, and show that when $c_0=0$ and $\delta>0$ the pressure adjoint in general has no representative continuous in time.
\end{enumerate}

The paper is organized as follows. Section~\ref{sec:setting} introduces the weak formulation and the pressure-to-dilation operator. State well-posedness is considered in Section~\ref{sec:state}. Section~\ref{sec:ocp} formulates the optimal control problem. The adjoint systems and first-order optimality conditions are derived in Section~\ref{sec:adjoint}. The singular-limit estimates and the convergence of optimal controls are proved in Sections~\ref{sec:limits} and \ref{sec:optimal-limits}, respectively, and Section~\ref{sec:boundary-controls} treats boundary controls. Appendix~\ref{app:sharpness} collects the examples and proofs of the sharpness results.

\section{Model problem and functional setting}\label{sec:setting}
\subsection{Notation and weak formulation}
Let $T>0$ be a final time and let $\Omega\subset\mathbb R^d$, $d\in\{2,3\}$, be bounded, connected, and Lipschitz. Its boundary has a mechanical Dirichlet portion $\Gamma_{u,D}$ and a pressure Dirichlet portion $\Gamma_{p,D}$, each of positive surface measure; their complements $\Gamma_{u,N}:=\partial\Omega\setminus\Gamma_{u,D}$ and $\Gamma_{p,N}:=\partial\Omega\setminus\Gamma_{p,D}$ support traction and normal Darcy flux. We assume in addition that $\Gamma_{u,N}$ contains a nonempty relatively open subset of $\partial\Omega$; in particular, the skeleton is not fully clamped. Set
\begin{equation}
V=\{v\in H^1(\Omega)^d:\ v|_{\Gamma_{u,D}}=0\},\quad
Q=\{r\in H^1(\Omega):\ r|_{\Gamma_{p,D}}=0\},\quad
H=L^2(\Omega).
\label{eq:spaces}
\end{equation}
The norm $\norm\cdot_V$ is the $H^1(\Omega)^d$ norm, and $\norm\cdot_H$ and $(\cdot,\cdot)_H$ are the norm and inner product of $L^2(\Omega)$. Duality between $V^*,V$ and $Q^*,Q$ is indicated by $\dual\cdot\cdot$. The natural injection $Q\hookrightarrow H\simeq H^*\hookrightarrow Q^*$ is dense and continuous. For a Banach space $X$ we abbreviate $L^2X:=L^2(0,T;X)$, $L^\infty X:=L^\infty(0,T;X)$, and $H^1X:=H^1(0,T;X)$. Time derivatives are denoted by $\partial_t$ or by a subscript $t$, e.g.\ $u_t=\partial_tu$; for scalar functions of $t$ we also write $'$.

Throughout, $\bar c_0>0$ and $\bar\delta>0$ are fixed upper bounds for the storage and viscosity coefficients, and we consider $0\le c_0\le\bar c_0$ and $0\le\delta\le\bar\delta$.

Let $\mathcal C_{\rm e}$ and $\mathcal C_{\rm v}$ be time-independent, symmetric elasticity and viscosity tensors, uniformly bounded and positive definite on symmetric matrices, and let $\mathbf K=\mathbf K(x,t)$ be a symmetric permeability tensor with
\begin{equation}
\mathbf K\in W^{1,\infty}(0,T;L^\infty(\Omega)^{d\times d}),\qquad
\underline C_{\mathbf K} \abs\xi^2\le\mathbf K(x,t)\xi\cdot\xi\le \overline C_{\mathbf K}\abs\xi^2
\label{eq:K-assumption}
\end{equation}
for a.e.\ $(x,t)$ and all $\xi\in\mathbb R^d$, where $0<\underline C_{\mathbf K}\le \overline C_{\mathbf K}$. Define
\begin{align}
a_{\ee}(u,v)&=\int_\Omega\mathcal{C}_{\rm e} \eps(u):\eps(v)\dd x,&
a_{\vv}(u,v)&=\int_\Omega\mathcal{C}_{\rm v} \eps(u):\eps(v)\dd x,\nonumber\\
a_p(r,q)&=\int_\Omega\mathbf K\nabla r\cdot\nabla q\dd x,&
b(v,q)&=\alpha\int_\Omega q\divergence v\dd x,
\label{eq:forms}
\end{align}
where $\eps(v):=\frac12(\nabla v+\nabla v^\top)$ is the symmetric gradient and $\alpha>0$ is fixed. The form $a_p$ depends on $t$ through $\mathbf K$; we write $a_p(t;\cdot,\cdot)$ when the time must be shown and otherwise suppress it, and we set $a_p'(r,q):=\int_\Omega\partial_t\mathbf K\nabla r\cdot\nabla q\dd x$. We take $\norm q_Q:=a_p(0;q,q)^{1/2}$, the energy norm on $Q$ induced by the Darcy form. Korn's and Poincar\'e's inequalities and \eqref{eq:K-assumption} give constants independent of $c_0,\delta$ and $t$ such that
\begin{equation}
\begin{aligned}
	\underline C_{\ee}\norm v_V^2&\le a_{\ee}(v,v)\le\overline C_{\ee}\norm v_V^2, &\underline C_{\vv}\norm v_V^2&\le a_{\vv}(v,v)\le\overline C_{\vv}\norm v_V^2, 
	\\
	\underline C_p\norm q_Q^2&\le a_p(q,q)\le\overline C_p\norm q_Q^2, & \abs{a_p'(r,q)}&\le\kappa_1\norm r_Q\norm q_Q.
\end{aligned}
\label{eq:ellipticity}
\end{equation}
Let
\begin{equation}
A_{\ee},A_{\vv}:V\to V^*,\qquad \KK=\KK(t):Q\to Q^*,
\qquad B:V\to H,\quad Bv=\alpha\divergence v,
\label{eq:operator-def}
\end{equation}
be the operators associated with the forms $a_{\ee}$, $a_{\vv}$, $a_p(t;\cdot,\cdot)$, and $b$, that is, $\dual{A_{\ee}u}v=a_{\ee}(u,v)$, $\dual{A_{\vv}u}v=a_{\vv}(u,v)$, $\dual{\KK r}q=a_p(t;r,q)$, and $(Bv,q)_H=b(v,q)$. The adjoint $B^*:H\to V^*$ satisfies $\dual{B^*q}v=(Bv,q)_H$ for $q\in H$, and $H$ is identified with a subspace of $Q^*$ through the canonical embedding $H\hookrightarrow Q^*$; in this sense $B$ and $\KK^{-1}B$ also act from $V$ into $Q^*$ and $Q$, respectively. The operators $A_{\ee},A_{\vv},\KK$ are bounded isomorphisms between the indicated spaces. The permeability may depend on time as in \eqref{eq:K-assumption}; then $\KK=\KK(t)$ has bounds uniform in $t$, and $t\mapsto\KK(t)$ and $t\mapsto\KK(t)^{-1}$ are Lipschitz continuous in $\calL(Q,Q^*)$ and $\calL(Q^*,Q)$. We also use the $H$-realization $\KK_H(t)$ of $a_p(t;\cdot,\cdot)$, defined by
\begin{equation*}
D(\KK_H(t))=\bigl\{q\in Q:\ \exists f\in H\text{ with }a_p(t;q,r)=(f,r)_H\ \forall r\in Q\bigr\},\qquad
\KK_H(t)q=f.
\end{equation*}
It is positive and self-adjoint in $H$. Since $Q$ embeds compactly into $H$, $\KK_H(0)$ has a compact inverse, and its eigenfunctions form an orthonormal basis of $H$; moreover $\norm q_Q=a_p(0;q,q)^{1/2}$ is the form norm of $\KK_H(0)$. The derivative $a_p'$ enters only where a pressure equation is tested with a time derivative.

For weak mechanical traction $t_N$ and outward flux $\psi_N$, we use, for a.e.\ $t\in(0,T)$, the unified load notation
\begin{equation}
\dual{F(t)}v:=(f(t),v)_\Omega+(t_N(t),v)_{\Gamma_{u,N}},\qquad
\dual{G(t)}q:=(g(t),q)_\Omega-(\psi_N(t),q)_{\Gamma_{p,N}}.
\label{eq:loads}
\end{equation}
For example, for a.e.\ $t$, $f(t)\in L^2(\Omega)^d$ and $t_N(t)\in L^2(\Gamma_{u,N})^d$ define $F(t)\in V^*$ by the trace theorem. Analogously, $g(t)\in H$ and $\psi_N(t)\in L^2(\Gamma_{p,N})$ define $G(t)\in Q^*$. With outward unit normal $\mathbf n$, these loads correspond to the boundary conditions
\begin{equation*}
\begin{aligned}
u(t)&=0\ \text{on }\Gamma_{u,D}, &\qquad
\bigl(\mathcal C_{\rm e}\eps(u(t))+\delta\mathcal C_{\rm v}\eps(u_t(t))-\alpha p(t)I\bigr)\mathbf n&=t_N(t)\ \text{on }\Gamma_{u,N},\\
p(t)&=0\ \text{on }\Gamma_{p,D}, &\qquad
-\mathbf K(t)\nabla p(t)\cdot\mathbf n&=\psi_N(t)\ \text{on }\Gamma_{p,N},
\end{aligned}
\end{equation*}
where the natural conditions on $\Gamma_{u,N}$ and $\Gamma_{p,N}$ hold in the variational sense.

\subsection{The pressure-to-dilation operator}\label{subsec:mass}
The operator
\begin{equation}
M_0:=BA_{\ee}^{-1}B^*:H\to H,\qquad
M_{c_0}:=c_0I+M_0
\label{eq:mass}
\end{equation}
is bounded, self-adjoint and nonnegative. Indeed,
\begin{equation}
(M_0q,q)_H=a_{\ee}(A_{\ee}^{-1}B^*q,A_{\ee}^{-1}B^*q)
\quad(q\in H).
\label{eq:mass-pos}
\end{equation}
For estimates that are uniform as $c_0\downarrow0$, we use the following pressure-to-dilation condition, which holds under the standing assumption on $\Gamma_{u,N}$.
\begin{lemma}[full-pressure inf--sup]\label{lem:infsup}
There is $\beta>0$ such that
\begin{equation}
\inf_{q\in H\setminus\{0\}}\sup_{v\in V\setminus\{0\}}
\frac{(Bv,q)_H}{\norm v_V\norm q_H}\ge\beta.
\label{eq:infsup}
\end{equation}
\end{lemma}
\begin{proof}
By \cite[Lemma~B.1]{BLMW2017}, every $q\in H$ admits $v\in V$ with $Bv=q$ and $\norm v_V\le C\norm q_H$; hence \eqref{eq:infsup} holds with $\beta=C^{-1}$.
\end{proof}
\begin{lemma}[uniform coercivity of $M_{c_0}$]\label{lem:mass-coercivity}
Lemma~\ref{lem:infsup} implies
\begin{equation}
(M_{c_0}q,q)_H\ge c_*\norm q_H^2,
\qquad c_*:=\frac{\beta^2}{\overline C_{\ee}}>0,
\qquad c_0\ge0.
\label{eq:uniform-mass}
\end{equation}
In particular all $M_{c_0}$ are uniformly bounded positive isomorphisms on $H$ for $0\le c_0\le\bar c_0$.
\end{lemma}
\begin{proof}
By the dual norm characterized by $A_{\ee}$,
\[
(M_0q,q)_H
=\sup_{v\in V\setminus\{0\}}\frac{\abs{(Bv,q)_H}^2}{a_{\ee}(v,v)}
\ge\frac1{\overline C_{\ee}}
\sup_{v\ne0}\frac{\abs{(Bv,q)_H}^2}{\norm v_V^2}
\ge c_*\norm q_H^2.
\]
Adding $c_0\norm q_H^2$ proves the first assertion. The upper bounds for $M_{c_0}$ on $0\le c_0\le\bar c_0$ follow directly from boundedness of $B$ and $A_{\ee}^{-1}$.
\end{proof}

The letter $C$ denotes a generic positive constant, possibly different at each occurrence. Unless stated otherwise, it may depend on $\Omega$ and its boundary partition, $T$, $\alpha$, the constants in \eqref{eq:ellipticity} (including $\kappa_1$), $\beta$ and hence $c_*$, $\norm B$, $\bar c_0$, and $\bar\delta$, but not on $c_0$ or $\delta$; a constant that may depend on a parameter carries it as a subscript, as in $C_\delta$.

\subsection{The four parameter regimes}
For $c_0,\delta\ge0$, the weak state system is, for a.e.\ $t\in(0,T)$,
\begin{align}
\delta a_{\vv}(u_t,v)+a_{\ee}(u,v)-(p,Bv)_H&=\dual Fv,\qquad v\in V,
\label{eq:state-mech}\\
\dual{\partial_t(c_0p+Bu)}q+a_p(t;p,q)&=\dual Gq,\qquad q\in Q,
\label{eq:state-fluid}
\end{align}
For $\delta=0$ the viscous terms are absent, and no time derivative of $u$ is required.
In operator notation,
\begin{equation}
\delta A_{\vv}u_t+A_{\ee}u-B^*p=F,
\qquad (c_0p+Bu)_t+\KK p=G.
\label{eq:state-ops}
\end{equation}
The independent initial data depend on the parameter regime. The natural choices are summarized below.

\begin{center}
\renewcommand{\arraystretch}{1.22}
\begin{tabular}{@{}lll@{}}
\toprule
Regime & Natural evolutionary data & Evolution operator\\
\midrule
$c_0>0,\ \delta>0$ & $u(0)=u_0,\ p(0)=p_0$ & $\bigl(\begin{smallmatrix}\delta A_{\vv}&0\\ B&c_0I\end{smallmatrix}\bigr)$ acting on $(u,p)$\\
$c_0=0,\ \delta>0$ & $u(0)=u_0$ & $\delta A_{\vv}+B^*\KK^{-1}B$\\
$c_0>0,\ \delta=0$ & $p(0)=p_0$ or fluid content & $M_{c_0}=c_0I+BA_{\ee}^{-1}B^*$\\
$c_0=0,\ \delta=0$ & $p(0)=p_0$ or fluid content & $M_0=BA_{\ee}^{-1}B^*$\\
\bottomrule
\end{tabular}
\end{center}
At $\delta=0$ the initial displacement is not independent: if $p(0)=p_0$ is prescribed, it is necessarily
\begin{equation}
u(0)=A_{\ee}^{-1}\bigl(F(0)+B^*p_0\bigr).
\label{eq:elastic-compatibility}
\end{equation}
If instead the initial fluid content $c_0p(0)+Bu(0)=\zeta_0\in H$ is prescribed, then \eqref{eq:elastic-compatibility} with $p_0=p(0)$ and Lemma~\ref{lem:mass-coercivity} give
\begin{equation*}
p(0)=M_{c_0}^{-1}\bigl(\zeta_0-BA_{\ee}^{-1}F(0)\bigr),
\end{equation*}
so the initial pressure then depends on $F(0)$. Unless stated otherwise, the initial pressure is prescribed. At $c_0=0$, $\delta>0$ a separate initial pressure is neither necessary nor prescribed.

\section{Well-posedness of the state equations}\label{sec:state}
We first establish the energy identity, then reduce the two singular regimes to standard Hilbert-space evolutions. The proofs use only energy-level spatial regularity.

\begin{lemma}[energy identity]\label{lem:state-energy}
Let $\delta>0$, and let $(u,p)$ satisfy \eqref{eq:state-ops} with $u\in H^1V$ and either $c_0>0$ and $p\in L^2Q\cap H^1Q^*$, or $c_0=0$ and $p\in L^2Q$. Then
\begin{equation}
\frac12\frac{\dd}{\dd t}\left\{a_{\ee}(u,u)+c_0\norm p_H^2\right\}
+\delta a_{\vv}(u_t,u_t)+a_p(p,p)
=\dual F{u_t}+\dual Gp
\label{eq:energy-identity}
\end{equation}
in the sense of distributions on $(0,T)$; equivalently, its integrated form holds on $[0,t]$ for every $t\in[0,T]$.
\end{lemma}
\begin{proof}
Test \eqref{eq:state-mech} by $v=u_t\in L^2V$ and \eqref{eq:state-fluid} by $q=p\in L^2Q$. The coupling terms are $-(p,Bu_t)_H$ and $(Bu_t,p)_H$, so they cancel. The symmetry and time independence of $a_{\ee}$ give $a_{\ee}(u,u_t)=\frac12\frac{\dd}{\dd t}a_{\ee}(u,u)$ for $u\in H^1V$. If $c_0>0$, the chain rule in the Gelfand triple $Q\hookrightarrow H\hookrightarrow Q^*$ gives $\dual{(c_0p)_t}p=\frac{c_0}2\frac{\dd}{\dd t}\norm p_H^2$; if $c_0=0$, this term is absent.
\end{proof}

\subsection{Positive storage and positive viscosity}
\begin{theorem}[State problem for $c_0>0$ and $\delta>0$]\label{thm:state-full}
Fix $c_0>0$ and $\delta>0$. Suppose
\begin{equation*}
F\in L^2V^*,\quad G\in L^2Q^*,\quad
u_0\in V,\quad p_0\in H.
\end{equation*}
There is a unique weak solution of \eqref{eq:state-ops} with initial data $(u_0,p_0)$ and
\begin{equation}
u\in H^1V,\qquad
p\in L^2Q\cap H^1Q^*\cap C([0,T];H).
\label{eq:full-state-space}
\end{equation}
It obeys
\begin{align}
&\norm u_{L^\infty V}^2+c_0\norm p_{L^\infty H}^2
 +\delta\norm{u_t}_{L^2V}^2+\norm p_{L^2Q}^2
\nonumber\\
&\qquad\le C\biggl(\norm{u_0}_V^2+c_0\norm{p_0}_H^2
 +\delta^{-1}\norm F_{L^2V^*}^2+\norm G_{L^2Q^*}^2\biggr),
\label{eq:state-full-est}
\end{align}
where $C$ depends only on the constants in \eqref{eq:ellipticity}, including $\overline C_{\ee}$ through the initial energy; in particular, it is independent of $c_0$, $\delta$, and the Galerkin discretization.
\end{theorem}
\begin{proof}
Choose nested finite-dimensional subspaces $V_n\subset V$ whose union is dense in $V$. For pressure take
$Q_n$ to be the span of the first $n$ eigenfunctions of the positive
self-adjoint reference Darcy realization $\KK_H(0)$ in $H$: the $H$-orthogonal
projection $\Pi_n:H\to Q_n$ is contractive in $H$ and in $\norm\cdot_Q$, which is the form norm
of $\KK_H(0)$. Choose $u_n(0)\in V_n$ such that $u_n(0)\to u_0$ in $V$, and set $p_n(0)=\Pi_np_0$, so that $p_n(0)\to p_0$ in $H$. The Galerkin equations are linear ordinary differential equations whose mass matrix is block lower triangular, with diagonal blocks $\delta a_{\vv}$ on $V_n$ and $c_0(\cdot,\cdot)_H$ on $Q_n$; for $c_0,\delta>0$ they have a unique absolutely continuous solution.

Testing the Galerkin equations by $\partial_tu_n\in V_n$ and $p_n\in Q_n$ as in the proof of Lemma~\ref{lem:state-energy} gives \eqref{eq:energy-identity} for $(u_n,p_n)$, and Young's inequality implies, with constants from \eqref{eq:ellipticity},
\begin{align*}
\abs{\dual F{\partial_t u_{n}}}
&\le \frac{\delta}{2}a_{\vv}(\partial_t u_{n},\partial_t u_{n})
   +\frac{C}{\delta}\norm F_{V^*}^2,\\
\abs{\dual G{p_n}}
&\le \frac12a_p(p_n,p_n)+C\norm G_{Q^*}^2.
\end{align*}
Integrating from $0$ to $t$, taking a supremum in $t$, and applying coercivity yields the uniform estimate \eqref{eq:state-full-est}. The Galerkin fluid equation reads
\begin{equation*}
c_0(\partial_t p_{n},q_n)_H=\dual G{q_n}-(B\partial_t u_{n},q_n)_H-a_p(t;p_n,q_n)
\qquad\text{for all }q_n\in Q_n,
\end{equation*}
For arbitrary $q\in Q$, take $q_n=\Pi_nq$ in this identity. Since
$(\partial_t p_{n},q)_H=(\partial_t p_{n},\Pi_n q)_H$ and
$\norm{\Pi_n q}_Q\le\norm q_Q$, this gives a bound of $\partial_t p_{n}$
in $L^2Q^*$ independent of $n$ for fixed $c_0>0$. Weak compactness gives a subsequence with $u_n\rightharpoonup u$ in $H^1V$ and $p_n\rightharpoonup p$ in $L^2Q\cap H^1Q^*$, and the limit satisfies the linear variational equations. Since $L^2Q\cap H^1Q^*\subset C([0,T];H)$ continuously \cite{LionsMagenes1972}, $p$ has the regularity \eqref{eq:full-state-space}. The maps $u\mapsto u(0)$ on $H^1V$ and $p\mapsto p(0)$ on $L^2Q\cap H^1Q^*$ are linear and continuous, hence weakly continuous; since $u_n(0)\to u_0$ in $V$ and $p_n(0)\to p_0$ in $H$, we obtain $u(0)=u_0$ and $p(0)=p_0$. Finally, Lemma~\ref{lem:state-energy} applied to the difference of two solutions, which has zero data, gives uniqueness.
\end{proof}
\begin{remark}
Estimate \eqref{eq:state-full-est} is not uniform as $\delta\downarrow0$. For the singular limit, we compare the viscous solution with an elastic solution in $H^1V$.
\end{remark}

\subsection{Zero storage and positive viscosity}
When $c_0=0$ the fluid equation \eqref{eq:state-fluid} determines
\begin{equation}
p=\KK^{-1}(G-Bu_t)\quad\text{in }Q.
\label{eq:elim-p-visc}
\end{equation}
Substitution into \eqref{eq:state-mech} gives the single displacement equation
\begin{equation}
N_\delta u_t+A_{\ee}u=F+B^*\KK^{-1}G,
\qquad N_\delta=N_\delta(t):=\delta A_{\vv}+B^*\KK(t)^{-1}B.
\label{eq:visc-reduced}
\end{equation}
Notice that $N_\delta(t):V\to V^*$ is symmetric and
\begin{equation}
\dual{N_\delta v}v
=\delta a_{\vv}(v,v)+\dual{Bv}{\KK^{-1}Bv}
\ge\delta\underline C_{\vv}\norm v_V^2,
\label{eq:N-coercive}
\end{equation}
where $Bv\in H\subset Q^*$ and $\KK^{-1}Bv\in Q$, so that $\dual{Bv}{\KK^{-1}Bv}\ge0$. By \eqref{eq:K-assumption}, $t\mapsto N_\delta(t)$ and $t\mapsto N_\delta(t)^{-1}$ are Lipschitz continuous in operator norm.

\begin{theorem}[State problem for $c_0=0$ and $\delta>0$]\label{thm:state-visc}
Fix $\delta>0$ and $c_0=0$. For $F\in L^2V^*$, $G\in L^2Q^*$, and $u_0\in V$, there exists a unique pair
\begin{equation}
u\in H^1V,\qquad p\in L^2Q,
\label{eq:visc-state-space}
\end{equation}
satisfying \eqref{eq:state-ops} and $u(0)=u_0$. The pressure is given by \eqref{eq:elim-p-visc}, and
\begin{equation}
\norm u_{L^\infty V}^2+\delta\norm{u_t}_{L^2V}^2+\norm p_{L^2Q}^2
\le C_\delta\bigl(\norm{u_0}_V^2+\norm F_{L^2V^*}^2+\norm G_{L^2Q^*}^2\bigr).
\label{eq:visc-est}
\end{equation}
No pressure initial trace is prescribed.
\end{theorem}
\begin{proof}
By \eqref{eq:N-coercive}, $N_\delta$ is a bounded symmetric isomorphism $V\to V^*$. Equation \eqref{eq:visc-reduced} is an ordinary evolution equation on the Hilbert space $V$ after applying $N_\delta^{-1}$:
\begin{equation*}
u_t+N_\delta^{-1}A_{\ee}u=N_\delta^{-1}(F+B^*\KK^{-1}G).
\end{equation*}
The operator $N_\delta^{-1}A_{\ee}$ is bounded on $V$ uniformly in $t$ and continuous in $t$ by \eqref{eq:K-assumption}, and the right-hand side belongs to $L^2V$; so the linear Hilbert-space ODE theorem yields a unique $u\in H^1V$ for each $u_0\in V$. Definition \eqref{eq:elim-p-visc} gives $p\in L^2Q$ and verifies both original equations. Lemma~\ref{lem:state-energy} and Young's inequality, as in the proof of Theorem~\ref{thm:state-full}, give \eqref{eq:visc-est}. The same argument for homogeneous data proves uniqueness.
\end{proof}

\subsection{Zero viscosity}
At $\delta=0$ the mechanical equilibrium equation is instantaneous:
\begin{equation}
u=A_{\ee}^{-1}(F+B^*p).
\label{eq:u-from-p}
\end{equation}
Since $u=A_{\ee}^{-1}(F+B^*p)$, the fluid content is $c_0p+Bu=M_{c_0} p+BA_{\ee}^{-1}F$, so pressure satisfies
\begin{equation}
(M_{c_0} p)_t+\KK p=G-BA_{\ee}^{-1}F_t\quad\text{in }Q^*.
\label{eq:pressure-reduced}
\end{equation}

\begin{theorem}[State problem for $c_0\ge0$ and $\delta=0$]\label{thm:state-storage}
Let $0\le c_0\le\bar c_0$, $\delta=0$, $F\in H^1V^*$, $G\in L^2Q^*$, and $p_0\in H$. Then there is a unique solution of \eqref{eq:pressure-reduced} with
\begin{equation}
p\in L^2Q\cap C([0,T];H),\qquad
(M_{c_0} p)_t\in L^2Q^*,
\label{eq:pressure-state-space}
\end{equation}
$p(0)=p_0$, and $u=A_{\ee}^{-1}(F+B^*p)\in C([0,T];V)$. Moreover,
\begin{align}
\sup_{0\le t\le T}(M_{c_0} p(t),p(t))_H+\norm p_{L^2Q}^2
\le C\left((M_{c_0} p_0,p_0)_H+\norm G_{L^2Q^*}^2
+\norm{F_t}_{L^2V^*}^2\right),
\label{eq:pressure-state-est}
\end{align}
where $C$ depends only on $\norm B$, $\underline C_{\ee}$, $\underline C_p$, and the embedding constant of $Q\hookrightarrow H$; in particular, it is independent of $c_0$. Moreover $\norm p_{L^\infty H}^2\le c_*^{-1}\sup_t(M_{c_0}p(t),p(t))_H$.
\end{theorem}
\begin{proof}
By Lemma~\ref{lem:mass-coercivity}, $M_{c_0}$ is a bounded positive isomorphism on $H$ for every $c_0\ge0$, including $c_0=0$. It equips $H$ with the equivalent inner product $(x,y)_{M_{c_0}}=(M_{c_0} x,y)_H$. Consider the triple $Q\hookrightarrow H_{M_{c_0}}\hookrightarrow Q^*$. With this inner product, the pivot embedding $H_{M_{c_0}}\hookrightarrow Q^*$ is $x\mapsto(M_{c_0}x,\cdot)_H$. Since $M_{c_0}$ is an isomorphism on $H$, this embedding has the same range as the standard embedding $H\hookrightarrow Q^*$, and the time derivative of $p$ in the weighted triple is $(M_{c_0}p)_t$. The bilinear form $a_p(t;\cdot,\cdot)$ is measurable in $t$ and continuous and coercive on $Q$ uniformly in $t$, as allowed in the parabolic existence theorem. The right-hand side $G_F:=G-BA_{\ee}^{-1}F_t$ belongs to $L^2Q^*$, since
\begin{equation*}
\norm{BA_{\ee}^{-1}F_t}_{L^2H}\le \norm B\norm{A_{\ee}^{-1}}\norm{F_t}_{L^2V^*}.
\end{equation*}
The standard Hilbert-triple parabolic existence theorem (see, e.g., \cite{LionsMagenes1972}), which can also be proved by a Galerkin scheme based on $Q$, now gives a unique
$p\in L^2Q\cap C([0,T];H_{M_{c_0}})$ with $(M_{c_0} p)_t\in L^2Q^*$ and $p(0)=p_0$.

Testing the equation with $p\in L^2Q$ and using the chain rule in the triple gives
\begin{equation*}
\frac12\frac{\dd}{\dd t}(M_{c_0} p,p)_H+a_p(t;p,p)=\dual{G_F}p.
\end{equation*}
Young's inequality in the form $\abs{\dual{G_F}p}\le\frac12a_p(t;p,p)+\frac1{2\underline C_p}\norm{G_F}_{Q^*}^2$, integration, and the above bound on $BA_{\ee}^{-1}F_t$ yield \eqref{eq:pressure-state-est}; no constant in this argument depends on $c_0$. Finally \eqref{eq:u-from-p} and continuity of $F$ in $V^*$ yield the claimed regularity of $u$.
\end{proof}
\begin{remark}[Pressure time regularity]
Theorem~\ref{thm:state-storage} directly controls $(M_{c_0} p)_t$ in $Q^*$. Bounded invertibility of $M_{c_0}$ on $H$ does \emph{not} alone imply a bounded realization $M_{c_0}^{-1}:Q^*\to Q^*$. A stronger statement $p_t\in L^2H$ is established in Lemma~\ref{lem:uniform-reg} when the right-hand side of \eqref{eq:pressure-reduced} belongs to $L^2H$ and $p_0\in Q$.
\end{remark}

\section{The optimal control problem}\label{sec:ocp}
We use the standard Hilbert-space framework for linear-quadratic optimal control problems; see \cite{Lions1971,Troeltzsch2010}. Let $\UU$ be a real Hilbert control space and $\UU_{\ad}\subset\UU$ a nonempty, closed, convex subset. For fixed $(c_0,\delta)$ let
\begin{equation*}
\mathcal F_\delta:=
\begin{cases}
L^2V^*,&\delta>0,\\
H^1V^*,&\delta=0,
\end{cases}
\qquad
\mathcal G:=L^2Q^*
\end{equation*}
be the load spaces of Theorems~\ref{thm:state-full}, \ref{thm:state-visc}, and \ref{thm:state-storage}. The loads in \eqref{eq:loads} are taken to be affine in the control,
\begin{equation}
F=F_b+C_Fm,\qquad G=G_b+C_Gm,\qquad m\in\UU,
\label{eq:affine-data}
\end{equation}
where the background loads $F_b\in\mathcal F_\delta$ and $G_b\in\mathcal G$ do not depend on the control, and the control operators $C_F\in\calL(\UU,\mathcal F_\delta)$ and $C_G\in\calL(\UU,\mathcal G)$ are linear and bounded. The initial data required by the relevant state theorem are fixed independently of the control unless explicitly stated otherwise. Let
\begin{equation}
\YY:=L^2(0,T;L^2(\Omega)^d)\times L^2H,
\qquad
\YY_V:=L^2V\times L^2H.
\label{eq:tracking-spaces}
\end{equation}
Every state obtained in Section~\ref{sec:state} belongs to $\YY$. Let $\mathcal S_{c_0,\delta}m$ denote the state generated by the loads $(C_Fm,C_Gm)$ with homogeneous initial data, that is, with the independent initial quantities of the regime set to zero ($u(0)$ and $p(0)$ for $c_0>0,\delta>0$; $u(0)$ for $c_0=0,\delta>0$; and the prescribed initial pressure or initial fluid content for $\delta=0$), and let $y^b_{c_0,\delta}$ denote the background state generated by $(F_b,G_b)$ and the fixed initial data. By linearity of \eqref{eq:state-ops} and the a priori estimates of Theorems~\ref{thm:state-full}, \ref{thm:state-visc}, and \ref{thm:state-storage},
\begin{equation}
(u(m),p(m))=\mathcal S_{c_0,\delta}m+y^b_{c_0,\delta},\qquad
\mathcal S_{c_0,\delta}\in\calL(\UU,\YY).
\label{eq:control-map}
\end{equation}

For a desired trajectory $y_d=(u_d,p_d)\in\YY$ and a fixed $\gamma>0$, consider
\begin{equation}
\min_{m\in\UU_{\ad}}j_{c_0,\delta}(m),\quad
j_{c_0,\delta}(m)=\frac12\norm{\mathcal S_{c_0,\delta}m+y^b_{c_0,\delta}-y_d}_{\YY}^2
+\frac\gamma2\norm m_\UU^2.
\label{eq:cost}
\end{equation}
We use displacement and pressure tracking in $\YY$; the same argument applies to an equivalent weighted inner product. Throughout, $\mathcal S_{c_0,\delta}^*\in\calL(\YY,\UU)$ denotes the Hilbert adjoint of $\mathcal S_{c_0,\delta}$ with respect to the inner product of $\YY$; the smaller space $\YY_V$ is used only for estimates.

\begin{theorem}[Existence and uniqueness of an optimal control]\label{thm:ocp-unique}
In the setting of \eqref{eq:affine-data}--\eqref{eq:control-map}, problem \eqref{eq:cost} has a unique minimizer $\bar m_{c_0,\delta}\in\UU_{\ad}$. It is characterized by
\begin{equation}
\bigl(\mathcal S_{c_0,\delta}^*(\mathcal S_{c_0,\delta}\bar m_{c_0,\delta}+y^b_{c_0,\delta}-y_d)
+\gamma\bar m_{c_0,\delta},\ m-\bar m_{c_0,\delta}\bigr)_\UU\ge0
\quad(m\in\UU_{\ad}).
\label{eq:abstract-VI}
\end{equation}
\end{theorem}
\begin{proof}
The reduced objective is continuous and quadratic on $\UU$, with
\begin{equation*}
j_{c_0,\delta}''(m)[r,r]
=\norm{\mathcal S_{c_0,\delta}r}_{\YY}^2+\gamma\norm r_\UU^2
\ge\gamma\norm r_\UU^2,\qquad r\in\UU.
\end{equation*}
Since $j_{c_0,\delta}(m)\ge\frac\gamma2\norm m_\UU^2$, every minimizing sequence in $\UU_{\ad}$ is bounded and has a weakly convergent subsequence. The set $\UU_{\ad}$ is closed and convex, hence weakly closed, and $j_{c_0,\delta}$ is continuous and convex, hence weakly lower semicontinuous; therefore a minimizer exists. The Hessian bound above shows that $j_{c_0,\delta}$ is $\gamma$-strongly convex, so the minimizer is unique. Its derivative is $\mathcal S_{c_0,\delta}^*(\mathcal S_{c_0,\delta}m+y^b_{c_0,\delta}-y_d)+\gamma m$, and \eqref{eq:abstract-VI} is the first-order condition, which is necessary and, by convexity, sufficient; see, e.g., \cite[Chapter~2]{Troeltzsch2010}.
\end{proof}

\begin{remark}[Control spaces]\label{rem:control-spaces}
If $\delta>0$, mechanical forcing and traction controls may be taken in $L^2V^*$, or in a Hilbert space mapped continuously into it. For $\delta=0$, the state theory used here requires $F\in H^1V^*$. Fluid source and flux controls remain admissible in $L^2Q^*$ in all four regimes. Section~\ref{sec:optimal-limits} therefore considers distributed fluid-source controls in a common $L^2$ control space.
\end{remark}

\section{Adjoint equations and optimality conditions}\label{sec:adjoint}
\subsection{Adjoint formulation}\label{subsec:adj-form}
Let the observation residuals be
\begin{equation}
R_u\in L^2V^*,\qquad R_p\in L^2H.
\label{eq:adj-data}
\end{equation}
In the optimal control problem \eqref{eq:cost}, $R_u=u-u_d\in L^2(0,T;L^2(\Omega)^d)\hookrightarrow L^2V^*$ and $R_p=p-p_d\in L^2H$.
Let $\tilde u$ denote the mechanical adjoint and $\tilde p$ the pressure adjoint. Formal integration by parts in time in the \emph{original} forward system \eqref{eq:state-ops} leads to the backward system
\begin{align}
-\bigl(\delta A_{\vv}\tilde u+B^*\tilde p\bigr)_t+A_{\ee}\tilde u&=R_u\quad\text{in }V^*,
\label{eq:adj-mech}\\
-\bigl(c_0\tilde p\bigr)_t-B\tilde u+\KK\tilde p&=R_p\quad\text{in }Q^*,
\label{eq:adj-fluid}
\end{align}
understood in the sense of distributions on $(0,T)$. Only the combined quantities $\delta A_{\vv}\tilde u+B^*\tilde p$ and $c_0\tilde p$ are differentiated in time; time derivatives of $\tilde u$ and $\tilde p$ separately are used only in regimes in which the corresponding regularity is established (Theorems~\ref{thm:adj-full} and \ref{thm:adj-elastic}). For sufficiently smooth functions, pairing \eqref{eq:adj-mech}--\eqref{eq:adj-fluid} with a forward variation $(w,r)$ produces the temporal endpoint contribution
\begin{equation}
\left[\dual{\delta A_{\vv}\tilde u+B^*\tilde p}{w}+(c_0\tilde p,r)_H\right]_{t=0}^{t=T}
=\left[\delta a_{\vv}(w,\tilde u)+(Bw+c_0r,\tilde p)_H\right]_{t=0}^{t=T}.
\label{eq:adj-terminal-form}
\end{equation}
We impose the terminal conditions
\begin{equation}
\bigl(\delta A_{\vv}\tilde u+B^*\tilde p\bigr)(T)=0\quad\text{in }V^*,
\qquad \bigl(c_0\tilde p\bigr)(T)=0\quad\text{in }H,
\label{eq:terminal-general}
\end{equation}
which make the contribution at $t=T$ vanish for every forward variation. Throughout, we write
\begin{equation*}
z:=\delta A_{\vv}\tilde u+B^*\tilde p,
\end{equation*}
so that \eqref{eq:terminal-general} reads $z(T)=0$ and $(c_0\tilde p)(T)=0$. The contribution at $t=0$ is determined by the initial data of the variation; see Proposition~\ref{prop:duality}. In a degenerate regime only a combined quantity may possess a terminal trace. The table below summarizes the independent state data and the resulting adjoint terminal conditions; the last row uses Lemma~\ref{lem:infsup}.
\begin{center}
\renewcommand{\arraystretch}{1.22}
\begin{tabular}{@{}lll@{}}
\toprule
Regime & State initial data & Adjoint terminal data\\
\midrule
$c_0>0,\ \delta>0$ & $u(0),\ p(0)$ & $\tilde u(T)=0,\ \tilde p(T)=0$\\
$c_0=0,\ \delta>0$ & $u(0)$ & $(\delta A_{\vv}\tilde u+B^*\tilde p)(T)=0$\\
$c_0>0,\ \delta=0$ & $p(0)$ or $c_0p(0)+Bu(0)$ & $\tilde p(T)=0$\\
$c_0=0,\ \delta=0$ & $p(0)$ or $Bu(0)$ & $\tilde p(T)=0$\\
\bottomrule
\end{tabular}
\end{center}
The regimes $c_0>0,\delta>0$, $\delta=0$, and $c_0=0,\delta>0$ are treated separately below.

\subsection{The $c_0>0,\delta>0$ adjoint}
\begin{theorem}[Adjoint problem for $c_0>0$ and $\delta>0$]\label{thm:adj-full}
Fix $c_0,\delta>0$. For data \eqref{eq:adj-data}, equations \eqref{eq:adj-mech}--\eqref{eq:adj-fluid} with \eqref{eq:terminal-general} have a unique solution
\begin{equation}
\tilde u\in H^1V,\qquad
\tilde p\in H^1H\cap L^\infty Q,\qquad
\KK_H\tilde p\in L^2H.
\label{eq:full-adj-space}
\end{equation}
There is a constant $C_{c_0,\delta}$ such that
\begin{align}
&\norm {\tilde u}_{L^2V}^2+\delta\norm {\tilde u}_{L^\infty V}^2
+c_0\norm{\tilde p_t}_{L^2H}^2+\norm {\tilde p}_{L^\infty Q}^2
\le C_{c_0,\delta}\bigl(\norm{R_u}_{L^2V^*}^2+\norm{R_p}_{L^2H}^2\bigr).
\label{eq:full-adj-est}
\end{align}
\end{theorem}
\begin{proof}
Use $V_n$ and the Darcy eigenfunction spaces $Q_n$ as in the proof of
Theorem~\ref{thm:state-full}; reverse time and solve from the zero
initial coefficient vectors. After time reversal the block time-derivative matrix is block upper triangular,
\[
\begin{pmatrix}\delta A_{\vv,n}&B_n^*\\0&c_0I_n\end{pmatrix},
\]
where $B_n^*:Q_n\to V_n^*$ is defined by $\dual{B_n^*q_n}{v_n}=(q_n,Bv_n)_H$ for $q_n\in Q_n$ and $v_n\in V_n$. Its diagonal blocks are coercive, so it is invertible. It defines a unique finite-dimensional solution.

At the Galerkin level, test the Galerkin version of \eqref{eq:adj-mech} by $\tilde u_n$ and that of \eqref{eq:adj-fluid} by $-\partial_t\tilde p_n$. The terms $-(\partial_t\tilde p_n,B\tilde u_n)_H$ and $(B\tilde u_n,\partial_t\tilde p_n)_H$ cancel. We obtain
\begin{equation}
\begin{aligned}
&-\frac12\frac{\dd}{\dd t}\bigl\{\delta a_{\vv}(\tilde u_n,\tilde u_n)+a_p(t;\tilde p_n,\tilde p_n)\bigr\}
+a_{\ee}(\tilde u_n,\tilde u_n)+c_0\norm{\partial_t\tilde p_n}_H^2\\
&\qquad=\dual{R_u}{\tilde u_n}-(R_p,\partial_t\tilde p_n)_H-\frac12a_p'(\tilde p_n,\tilde p_n).
\end{aligned}
\label{eq:adj-full-identity}
\end{equation}
Integrate over $(t,T)$ and use $\tilde u_n(T)=0$ and $\tilde p_n(T)=0$. By \eqref{eq:ellipticity}, $\abs{a_p'(\tilde p_n,\tilde p_n)}\le\kappa_1\underline C_p^{-1}a_p(\tilde p_n,\tilde p_n)$. Apply the backward Gronwall inequality and use
\begin{align*}
\abs{\dual{R_u}{\tilde u_n}}&\le\tfrac12 a_{\ee}(\tilde u_n,\tilde u_n)+C\norm{R_u}_{V^*}^2,\\
\abs{(R_p,\partial_t\tilde p_n)_H}&\le\tfrac {c_0}2\norm{\partial_t\tilde p_n}_H^2+C c_0^{-1}\norm{R_p}_H^2.
\end{align*}
This proves \eqref{eq:full-adj-est} for $(\tilde u_n,\tilde p_n)$, uniformly in $n$. Testing the Galerkin version of \eqref{eq:adj-mech} by $\partial_t\tilde u_n$ gives moreover
\begin{equation*}
\delta\underline C_{\vv}\norm{\partial_t\tilde u_n}_V
\le\overline C_{\ee}\norm{\tilde u_n}_V+\norm B\norm{\partial_t\tilde p_n}_H+\norm{R_u}_{V^*},
\end{equation*}
so $\partial_t\tilde u_n$ is bounded in $L^2V$ for fixed $c_0,\delta>0$.

By weak and weak-$*$ compactness, a subsequence satisfies $\tilde u_n\rightharpoonup\tilde u$ in $H^1V$, and $\tilde p_n\rightharpoonup\tilde p$ in $H^1H$ and weakly-$*$ in $L^\infty Q$. The limit satisfies
\begin{equation*}
-\delta A_{\vv}\tilde u_t-B^*\tilde p_t+A_{\ee}\tilde u=R_u\ \text{in }L^2V^*,\qquad
-c_0\tilde p_t-B\tilde u+\KK\tilde p=R_p\ \text{in }L^2Q^*,
\end{equation*}
which is \eqref{eq:adj-mech}--\eqref{eq:adj-fluid} for this regularity. The maps $v\mapsto v(T)$ on $H^1V$ and $q\mapsto q(T)$ on $H^1H$ are linear and continuous, hence weakly continuous; since $\tilde u_n(T)=0$ and $\tilde p_n(T)=0$, we obtain $\tilde u(T)=0$ and $\tilde p(T)=0$, which implies \eqref{eq:terminal-general}. The bound \eqref{eq:full-adj-est} passes to the limit by weak and weak-$*$ lower semicontinuity. By \eqref{eq:adj-fluid}, $\KK\tilde p=R_p+B\tilde u+c_0\tilde p_t$ in $L^2Q^*$, and the right-hand side lies in $L^2H$; hence, by the definition of $\KK_H(t)$ in Section~\ref{sec:setting}, $\tilde p(t)\in D(\KK_H(t))$ for a.e.\ $t$ and $\KK_H\tilde p\in L^2H$. This proves existence with the regularity \eqref{eq:full-adj-space}.

For uniqueness let $R_u=R_p=0$, let $(\tilde u,\tilde p)$ be a solution with the regularity \eqref{eq:full-adj-space}, and set
\begin{equation*}
h(t):=\int_t^T\tilde u(s)\dd s,\qquad\text{so that}\qquad h\in H^1V,\quad h_t=-\tilde u,\quad h(T)=0.
\end{equation*}
Since $c_0>0$, \eqref{eq:terminal-general} gives $z(T)=0$ and $\tilde p(T)=0$. By \eqref{eq:adj-mech}, $z_t=A_{\ee}\tilde u$ in $L^2V^*$, and integration over $(t,T)$ gives $z(t)=-A_{\ee}h(t)$. Substituting $z=\delta A_{\vv}\tilde u+B^*\tilde p$ and $\tilde u=-h_t$ into $z=-A_{\ee}h$, and $\tilde u=-h_t$ into \eqref{eq:adj-fluid}, we obtain
\begin{equation*}
-\delta A_{\vv}h_t+A_{\ee}h+B^*\tilde p=0\ \text{in }L^2V^*,\qquad
-c_0\tilde p_t+Bh_t+\KK\tilde p=0\ \text{in }L^2Q^*,
\end{equation*}
Test $-\delta A_{\vv}h_t+A_{\ee}h+B^*\tilde p=0$ by $-h_t\in L^2V$ and $-c_0\tilde p_t+Bh_t+\KK\tilde p=0$ by $\tilde p\in L^2Q$. The coupling terms $-(\tilde p,Bh_t)_H$ and $(Bh_t,\tilde p)_H$ cancel, and integration over $(t,T)$ with $h(T)=0$ and $\tilde p(T)=0$ gives
\begin{equation*}
\frac12a_{\ee}(h(t),h(t))+\frac{c_0}2\norm{\tilde p(t)}_H^2
+\int_t^T\Bigl\{\delta a_{\vv}(h_t,h_t)+a_p(s;\tilde p,\tilde p)\Bigr\}\dd s=0.
\end{equation*}
Hence $h=0$ and $\tilde p=0$, so $\tilde u=-h_t=0$.
\end{proof}

\subsection{The adjoint problem for $\delta=0$}
Let $\delta=0$ and $0\le c_0\le\bar c_0$. For $\tilde p\in H^1H$, \eqref{eq:adj-mech} gives
\begin{equation}
\tilde u=A_{\ee}^{-1}(R_u+B^*\tilde p_t),
\label{eq:adj-elim-m}
\end{equation}
and substitution into \eqref{eq:adj-fluid} yields
\begin{equation}
-M_{c_0}\tilde p_t+\KK \tilde p=R_p+BA_{\ee}^{-1}R_u=:\mathcal R,
\qquad \tilde p(T)=0.
\label{eq:adj-reduced-pressure}
\end{equation}
Here \eqref{eq:terminal-general} reduces to $\tilde p(T)=0$: for $c_0>0$ this is the condition $(c_0\tilde p)(T)=0$, and for $c_0=0$ the condition $z(T)=0$, which at $\delta=0$ reads $B^*\tilde p(T)=0$, implies $\tilde p(T)=0$ because $B^*$ is injective on $H$ by Lemma~\ref{lem:infsup}.
\begin{theorem}[Backward pressure equation]\label{thm:adj-elastic}
Let $R_u,R_p$ satisfy \eqref{eq:adj-data}, so that $\mathcal R\in L^2H$ in \eqref{eq:adj-reduced-pressure}. There exists a unique solution
\begin{equation}
\tilde p\in H^1H\cap L^\infty Q,\qquad \KK_H\tilde p\in L^2H,
\qquad \tilde u\in L^2V,
\label{eq:adj-elastic-space}
\end{equation}
where $\KK_H=\KK_H(t)$ is the $H$-realization of $a_p(t;\cdot,\cdot)$ from Section~\ref{sec:setting} and $\tilde u$ is defined by \eqref{eq:adj-elim-m}. If $0\le c_0\le\bar c_0$, then
\begin{equation}
\norm{\tilde p_t}_{L^2H}+\norm {\tilde p}_{L^\infty Q}+\norm{\KK_H\tilde p}_{L^2H}
+\norm {\tilde u}_{L^2V}
\le C\left(\norm{R_u}_{L^2V^*}+\norm{R_p}_{L^2H}\right),
\label{eq:adj-elastic-uniform}
\end{equation}
with $C$ depending on $T$ and $\kappa_1$ but independent of $c_0$.
\end{theorem}
\begin{proof}
Reverse time by writing $\widehat{\tilde p}(\tau)=\tilde p(T-\tau)$, $\widehat{\mathcal R}(\tau)=\mathcal R(T-\tau)$, and $\widehat\KK(\tau):=\KK(T-\tau)$. Then
\begin{equation*}
M_{c_0}\widehat{\tilde p}_\tau+\widehat\KK\widehat{\tilde p}=\widehat{\mathcal R},\qquad
\widehat{\tilde p}(0)=0.
\end{equation*}
By Lemma~\ref{lem:mass-coercivity}, $M_{c_0}$ is coercive on $H$ with the lower bound $c_*$, uniformly in $c_0\ge0$. Let $\widehat{\tilde p}_n$ be the Galerkin approximation of this equation in the spaces $Q_n$ of the proof of Theorem~\ref{thm:state-full}, with initial vector zero; it can be tested by $\partial_\tau\widehat{\tilde p}_n\in Q_n$. Let $\widehat a_p(\tau;\cdot,\cdot):=a_p(T-\tau;\cdot,\cdot)$, so that $\frac{\dd}{\dd\tau}\widehat a_p(\widehat{\tilde p}_n,\widehat{\tilde p}_n)=2\widehat a_p(\widehat{\tilde p}_n,\partial_\tau\widehat{\tilde p}_n)-a_p'(\widehat{\tilde p}_n,\widehat{\tilde p}_n)$, with $a_p'$ evaluated at $T-\tau$. Symmetry of $\widehat\KK(\tau)$ implies
\begin{equation}
(M_{c_0}\partial_\tau\widehat{\tilde p}_n,\partial_\tau\widehat{\tilde p}_n)_H
+\frac12\frac{\dd}{\dd\tau}\widehat a_p(\widehat{\tilde p}_n,\widehat{\tilde p}_n)
+\frac12a_p'(\widehat{\tilde p}_n,\widehat{\tilde p}_n)
=(\widehat{\mathcal R},\partial_\tau\widehat{\tilde p}_n)_H.
\label{eq:adj-elastic-diff}
\end{equation}
Apply the coercivity of $M_{c_0}$ and Young's inequality, bound the $a_p'$ term in absolute value by $\kappa_1\underline C_p^{-1}\widehat a_p(\widehat{\tilde p}_n,\widehat{\tilde p}_n)$, integrate from $0$ to $\tau$, use $\widehat{\tilde p}_n(0)=0$, and apply Gronwall's inequality to obtain
\begin{equation*}
\norm{\partial_\tau\widehat{\tilde p}_n}_{L^2H}^2
+\sup_{\tau\in[0,T]}\widehat a_p(\widehat{\tilde p}_n(\tau),\widehat{\tilde p}_n(\tau))
\le C\norm{\widehat{\mathcal R}}_{L^2H}^2.
\end{equation*}
The Galerkin solutions are therefore bounded in $H^1H$ and in $L^\infty Q$; a subsequence converges weakly in $H^1H$ and weakly-$*$ in $L^\infty Q$ to a limit that satisfies the equation in $L^2Q^*$, and $\widehat{\tilde p}(0)=0$ is preserved because $H^1H\subset C([0,T];H)$ and the Galerkin initial values vanish. The equation implies $\widehat\KK\widehat{\tilde p}=\widehat{\mathcal R}-M_{c_0}\widehat{\tilde p}_\tau\in L^2H$, and the estimate gives $\widehat{\tilde p}\in L^\infty Q$; moreover $\widehat{\tilde p}\in H^1H\subset C([0,T];H)$. Finally, \eqref{eq:adj-elim-m} yields $\norm {\tilde u}_{L^2V}\le C(\norm{R_u}_{L^2V^*}+\norm{\tilde p_t}_{L^2H})$. By Lemma~\ref{lem:mass-coercivity}, the upper and lower bounds for $M_{c_0}$ are uniform on $[0,\bar c_0]$; hence all constants are uniform in $c_0$. For uniqueness, let $(\tilde u_i,\tilde p_i)$, $i=1,2$, be two solutions with the same data $R_u,R_p$, and put $\tilde p:=\tilde p_1-\tilde p_2$. Then $\tilde p$ solves \eqref{eq:adj-reduced-pressure} with $\mathcal R=0$, and $\tilde u_1-\tilde u_2=A_{\ee}^{-1}B^*\tilde p_t$ by \eqref{eq:adj-elim-m}. Since $\tilde p\in H^1H\cap L^\infty Q$ and $M_{c_0}$ is symmetric and independent of time, testing $-M_{c_0}\tilde p_t+\KK\tilde p=0$ with $\tilde p$ gives $-\frac12\frac{\dd}{\dd t}(M_{c_0}\tilde p,\tilde p)_H+a_p(t;\tilde p,\tilde p)=0$. Integration over $(t,T)$ and $\tilde p(T)=0$ give $\frac12(M_{c_0}\tilde p(t),\tilde p(t))_H+\int_t^Ta_p(s;\tilde p,\tilde p)\dd s=0$, hence $\tilde p=0$ by Lemma~\ref{lem:mass-coercivity}, and therefore $\tilde u_1=\tilde u_2$.
\end{proof}
\begin{remark}[Effect of positive storage]
For $c_0>0$ the pressure mass is coercive with lower bound $c_0$ even without Lemma~\ref{lem:infsup}. Thus $R_p\in L^2H$ suffices to construct the adjoint and no time derivative of the pressure tracking residual is required. This is the fundamental simplifying feature of positive storage. The endpoint $c_0=\delta=0$ enjoys the same benefit through Lemma~\ref{lem:infsup}, that is, because the skeleton is not fully clamped; in particular $\tilde p(T)=0$ is then a legitimate terminal trace. Without the inf--sup condition, as for a fully clamped skeleton, the natural endpoint condition may specify only $B^*\tilde p(T)$, as in the more general degenerate theory \cite{Showalter1974,BS2022}.
\end{remark}

\subsection{The adjoint problem for $c_0=0$ and $\delta>0$}
Let $c_0=0$, $\delta>0$. Equations \eqref{eq:adj-mech}--\eqref{eq:adj-fluid} become
\begin{subequations}\label{eq:visc-adjoint-form}
\begin{align}
-\bigl(\delta A_{\vv}\tilde u+B^*\tilde p\bigr)_t+A_{\ee}\tilde u&=R_u,
\label{eq:visc-adjoint-mech}\\
-B\tilde u+\KK \tilde p&=R_p.
\label{eq:visc-adjoint-fluid}
\end{align}
\end{subequations}
With $z=\delta A_{\vv}\tilde u+B^*\tilde p$ as in Section~\ref{subsec:adj-form}, define
\begin{equation}
\ell:=B^*\KK^{-1}R_p\in L^2V^*.
\label{eq:z-ell}
\end{equation}
By \eqref{eq:visc-adjoint-fluid}, $\tilde p=\KK^{-1}(R_p+B\tilde u)$, so
\begin{equation}
z=N_\delta \tilde u+\ell,
\qquad \tilde u=N_\delta^{-1}(z-\ell).
\label{eq:recover-m}
\end{equation}

\begin{theorem}[Adjoint problem for $c_0=0$ and $\delta>0$]\label{thm:adj-visc}
Fix $\delta>0$, $c_0=0$, and assume \eqref{eq:adj-data}. There exists a unique adjoint pair
\begin{equation}
\tilde u\in L^2V,\qquad \tilde p\in L^2Q,
\label{eq:adj-visc-spaces}
\end{equation}
for which $z\in H^1V^*$, $z(T)=0$, and \eqref{eq:visc-adjoint-form} holds in the sense of distributions on $(0,T)$. Moreover,
\begin{equation}
\norm z_{H^1V^*}+\norm {\tilde u}_{L^2V}+\norm {\tilde p}_{L^2Q}
\le C_\delta\bigl(\norm{R_u}_{L^2V^*}+\norm{R_p}_{L^2H}\bigr).
\label{eq:adj-visc-est}
\end{equation}
No separate terminal values for $\tilde u$ or $\tilde p$ are asserted.
\end{theorem}
\begin{proof}
Substituting \eqref{eq:recover-m} into \eqref{eq:visc-adjoint-mech} gives the backward equation
\begin{equation}
-z_t+A_{\ee}N_\delta^{-1}z
=R_u+A_{\ee}N_\delta^{-1}\ell,\qquad z(T)=0.
\label{eq:z-equation}
\end{equation}
The right-hand side of \eqref{eq:z-equation} belongs to $L^2V^*$ because $R_u\in L^2V^*$, $R_p\in L^2H\hookrightarrow L^2Q^*$ and $\KK^{-1}:Q^*\to Q$ is bounded. Time reversal reduces \eqref{eq:z-equation} to a linear ODE on $V^*$ with zero initial data, whose operator $A_{\ee}N_\delta^{-1}$ is bounded uniformly in $t$ and continuous in $t$ by \eqref{eq:K-assumption}. The linear Hilbert-space ODE theorem and Gronwall's inequality give a unique $z\in H^1V^*$ and the corresponding estimate.

Equation \eqref{eq:recover-m} then gives $\tilde u\in L^2V$ and $\tilde p=\KK^{-1}(R_p+B\tilde u)\in L^2Q$, along with \eqref{eq:adj-visc-est}. Substitution verifies \eqref{eq:visc-adjoint-form}. In the homogeneous case $R_u=R_p=0$, the ODE has only $z=0$, and thus $\tilde u=\tilde p=0$. This proves uniqueness.
\end{proof}
\begin{theorem}[Pressure-adjoint time regularity at zero storage]\label{thm:adj-jump}
Let $c_0=0$, $\delta>0$, $R_u\in L^2V^*$, and
\begin{equation}
R_p=R_p^{\rm c}+\chi_{[t_*,T]}\,\xi,\qquad
R_p^{\rm c}\in C([0,T];H),\quad t_*\in(0,T),\quad \xi\in H\setminus\{0\}.
\label{eq:jump-data}
\end{equation}
Then the pressure adjoint $\tilde p$ of Theorem~\ref{thm:adj-visc} has no representative in $C([0,T];H)$. In particular, $\tilde p\notin H^1H$: positive structural viscosity alone does not ensure $H^1$ time regularity of the pressure adjoint for $L^2$ pressure tracking data when $c_0=0$.
\end{theorem}
\begin{proof}
By the definition of $z$, $\tilde u=\delta^{-1}A_{\vv}^{-1}(z-B^*\tilde p)$. Substituting this into \eqref{eq:visc-adjoint-fluid} gives
\begin{equation}
L_\delta\tilde p=R_p+\delta^{-1}BA_{\vv}^{-1}z,
\qquad L_\delta=L_\delta(t):=\KK(t)+\delta^{-1}BA_{\vv}^{-1}B^*:Q\to Q^*.
\label{eq:h-jump-decomp}
\end{equation}
For $q\in Q$,
\begin{equation*}
\dual{L_\delta q}q=a_p(t;q,q)+\delta^{-1}a_{\vv}\bigl(A_{\vv}^{-1}B^*q,A_{\vv}^{-1}B^*q\bigr)\ge\underline C_p\norm q_Q^2,
\end{equation*}
so $L_\delta(t)$ is a bounded, symmetric, coercive isomorphism $Q\to Q^*$, and by \eqref{eq:K-assumption} $t\mapsto L_\delta(t)^{-1}$ is continuous in $\calL(Q^*,Q)$. Hence
\begin{equation*}
\tilde p=\chi_{[t_*,T]}\,L_\delta^{-1}\xi
+L_\delta^{-1}\bigl(R_p^{\rm c}+\delta^{-1}BA_{\vv}^{-1}z\bigr).
\end{equation*}
Since $z\in H^1V^*\subset C([0,T];V^*)$ by Theorem~\ref{thm:adj-visc}, $BA_{\vv}^{-1}\in\calL(V^*,H)$, and $R_p^{\rm c}\in C([0,T];H)$, the second term belongs to $C([0,T];Q)\subset C([0,T];H)$. If $\tilde p$ had a representative in $C([0,T];H)$, so would $\chi_{[t_*,T]}L_\delta(\cdot)^{-1}\xi$; since $L_\delta(\cdot)^{-1}\xi\in C([0,T];Q)$, this forces $L_\delta(t_*)^{-1}\xi=0$, which is impossible because $\xi\ne0$. Thus $\tilde p$ has the nonzero jump $L_\delta(t_*)^{-1}\xi$ at $t_*$. For the scalar analogue $V=H=Q=\mathbb R$, $A_{\ee}=A_{\vv}=B=\KK=1$, one has $L_\delta=1+\delta^{-1}$, which recovers the explicit formula $\tilde p=(z+\delta R_p)/(1+\delta)$.
\end{proof}

\subsection{State--adjoint duality}
\begin{lemma}[Integration by parts]\label{lem:product}
Suppose $\Xi\in C([0,T];H)$ with $\Xi_t\in L^2Q^*$, and
$q\in H^1H\cap L^2Q$.
Then for every $0\le a\le b\le T$,
\begin{equation}
\int_a^b\bigl\{\dual{\Xi_t}{q}+(\Xi,q_t)_H\bigr\}\,dt
=(\Xi(b),q(b))_H-(\Xi(a),q(a))_H.
\label{eq:product-rule}
\end{equation}
\end{lemma}
\begin{proof}
Since $\Xi\in C([0,T];H)$ and $\Xi_t\in L^2Q^*$, we have $\Xi\in H^1Q^*\cap L^2H$. The Leibniz rule of \cite[Lemma~4.5]{ColliGilardiSprekels2018}, whose proof applies verbatim to the Hilbert triple $Q\hookrightarrow H\hookrightarrow Q^*$, shows that $t\mapsto(\Xi(t),q(t))_H$ is absolutely continuous with derivative $\dual{\Xi_t}q+(\Xi,q_t)_H$ a.e.; integration over $(a,b)$ gives \eqref{eq:product-rule}.
\end{proof}

We now state the state--adjoint identity at the level needed for optimization. For a control direction $\eta\in\UU$ we set $F_\eta:=C_F\eta$ and $G_\eta:=C_G\eta$, with $C_F,C_G$ from \eqref{eq:affine-data}. A \emph{control variation} is the forward solution $(w,r)$ with sources $(F_\eta,G_\eta)$ and homogeneous variations of the independent initial quantities: $w(0)=0$ and $r(0)=0$ for $c_0>0,\delta>0$; $w(0)=0$ for $c_0=0,\delta>0$; and, for $\delta=0$, either $r(0)=0$ or $c_0r(0)+Bw(0)=0$, according to whether the initial pressure or the initial fluid content is prescribed. In the first case for $\delta=0$, $w(0)=A_{\ee}^{-1}F_\eta(0)$ may be nonzero if the mechanical source changes at the initial time; in the second, $r(0)=-M_{c_0}^{-1}BA_{\ee}^{-1}F_\eta(0)$ may be nonzero.

\begin{proposition}[State--adjoint identity]\label{prop:duality}
Let $c_0,\delta\ge0$, and let $(w,r)$ be the solution, given by Theorem~\ref{thm:state-full}, \ref{thm:state-visc}, or \ref{thm:state-storage} according to the regime, of
\begin{equation}
\delta A_{\vv}w_t+A_{\ee}w-B^*r=F_w,\qquad \Xi_t+\KK r=G_w,\qquad \Xi:=c_0r+Bw,
\label{eq:variation-forward}
\end{equation}
for general loads $(F_w,G_w)$ in the data spaces of that theorem and for initial data admissible in that regime. Let $(\tilde u,\tilde p)$ be the adjoint solution of Theorem~\ref{thm:adj-full}, \ref{thm:adj-visc}, or \ref{thm:adj-elastic} with data \eqref{eq:adj-data}. Then
\begin{equation}
\begin{split}
\int_0^T\!\left\{\dual{R_u}w+(R_p,r)_H\right\}\dd t
&=\int_0^T\!\left\{\dual{F_w}{\tilde u}+\dual{G_w}{\tilde p}\right\}\dd t
+\mathcal E_0,
\end{split}
\label{eq:state-adj-duality}
\end{equation}
where the initial term is defined, in each regime, by quantities that exist in the corresponding solution spaces:
\begin{equation}
\mathcal E_0:=
\begin{cases}
\dual{z(0)}{w(0)}+c_0(r(0),\tilde p(0))_H, & c_0>0,\ \delta>0,\\
\dual{z(0)}{w(0)}, & c_0=0,\ \delta>0,\\
(\Xi(0),\tilde p(0))_H, & \delta=0.
\end{cases}
\label{eq:initial-term}
\end{equation}
For a control variation $(w,r)$ in the direction $\eta$, that is, for $(F_w,G_w)=(F_\eta,G_\eta)=(C_F\eta,C_G\eta)$ and homogeneous independent initial data, we write $\mathcal E_0(\eta)$. Then $\mathcal E_0(\eta)=0$ if $\delta>0$, or if $\delta=0$ and the initial fluid content $c_0p(0)+Bu(0)$ is fixed, as in \cite{BS2022}. If $\delta=0$ and the initial pressure is fixed, then
\begin{equation}
\mathcal E_0(\eta)
=\bigl(BA_{\ee}^{-1}F_\eta(0),\tilde p(0)\bigr)_H,
\label{eq:initial-correction}
\end{equation}
which vanishes for fluid-source and flux variations.
\end{proposition}
\begin{proof}
We argue directly with the weak solutions; $z$ is as in Section~\ref{subsec:adj-form}.

\emph{Case $\delta>0$.} Then $w\in H^1V$ by Theorems~\ref{thm:state-full} and \ref{thm:state-visc}, and $z\in H^1V^*$ with $z(T)=0$ by Theorems~\ref{thm:adj-full} and \ref{thm:adj-visc}. The product rule in $H^1V^*\times H^1V$ gives
\begin{equation*}
-\int_0^T\dual{z_t}w\dd t=\dual{z(0)}{w(0)}+\int_0^T\dual z{w_t}\dd t,
\qquad
\dual z{w_t}=\delta a_{\vv}(w_t,\tilde u)+(Bw_t,\tilde p)_H.
\end{equation*}
If $c_0>0$, then $c_0r\in C([0,T];H)$ with $(c_0r)_t\in L^2Q^*$, and $\tilde p\in H^1H\cap L^2Q$ with $\tilde p(T)=0$; Lemma~\ref{lem:product} gives
\begin{equation*}
-\int_0^T(c_0\tilde p_t,r)_H\dd t=c_0(r(0),\tilde p(0))_H+\int_0^T\dual{(c_0r)_t}{\tilde p}\dd t.
\end{equation*}
Pair \eqref{eq:adj-mech} with $w$ and \eqref{eq:adj-fluid} with $r$, insert the product-rule identity and, if $c_0>0$, the identity from Lemma~\ref{lem:product}, and collect terms:
\begin{align*}
\int_0^T\!\left\{\dual{R_u}w+(R_p,r)_H\right\}\dd t
&=\int_0^T\!\left\{\delta a_{\vv}(w_t,\tilde u)+a_{\ee}(w,\tilde u)-(r,B\tilde u)_H\right\}\dd t\\
&\quad+\int_0^T\!\left\{\dual{\Xi_t}{\tilde p}+a_p(r,\tilde p)\right\}\dd t
+\mathcal E_0,
\end{align*}
where for $c_0=0$ the term $(c_0\tilde p)_t$ is absent from \eqref{eq:adj-fluid}, Lemma~\ref{lem:product} is not used, and no trace of $r$ or $\tilde p$ occurs. The two integrals equal $\int_0^T\dual{F_w}{\tilde u}\dd t$ and $\int_0^T\dual{G_w}{\tilde p}\dd t$ since $\delta A_{\vv}w_t+A_{\ee}w-B^*r=F_w$ and $\Xi_t+\KK r=G_w$ by \eqref{eq:variation-forward}.

\emph{Case $\delta=0$.} Then $\tilde p\in H^1H\cap L^2Q$ with $\tilde p(T)=0$ by Theorem~\ref{thm:adj-elastic}, and \eqref{eq:adj-mech}--\eqref{eq:adj-fluid} read $A_{\ee}\tilde u-B^*\tilde p_t=R_u$ and $-c_0\tilde p_t-B\tilde u+\KK\tilde p=R_p$. By Theorem~\ref{thm:state-storage}, $r\in L^2Q\cap C([0,T];H)$ and $w=A_{\ee}^{-1}(F_w+B^*r)\in C([0,T];V)$, so $\Xi\in C([0,T];H)$ and $\Xi_t=G_w-\KK r\in L^2Q^*$. Pairing the adjoint equations with $w$ and $r$ and using $A_{\ee}w-B^*r=F_w$,
\begin{equation*}
\int_0^T\!\left\{\dual{R_u}w+(R_p,r)_H\right\}\dd t
=\int_0^T\!\left\{\dual{F_w}{\tilde u}-(\tilde p_t,\Xi)_H+a_p(r,\tilde p)\right\}\dd t.
\end{equation*}
By Lemma~\ref{lem:product}, $-\int_0^T(\tilde p_t,\Xi)_H\dd t=(\Xi(0),\tilde p(0))_H+\int_0^T\dual{\Xi_t}{\tilde p}\dd t$, and $\dual{\Xi_t}{\tilde p}+a_p(r,\tilde p)=\dual{G_w}{\tilde p}$. This is \eqref{eq:state-adj-duality} with $\delta=0$.

\emph{Initial terms.} Inserting the initial values of a control variation into \eqref{eq:initial-term} gives the stated values of $\mathcal E_0(\eta)$; for $\delta=0$ and fixed initial pressure, $\Xi(0)=Bw(0)=BA_{\ee}^{-1}F_\eta(0)$.
\end{proof}

\begin{theorem}[First-order optimality conditions]\label{thm:optimality}
Let $\bar m$ solve \eqref{eq:cost}, let $(\bar u,\bar p)$ be its state, and solve the corresponding adjoint with $R_u=\bar u-u_d$ and $R_p=\bar p-p_d$. Then
\begin{equation}
\int_0^T\!\left\{\dual{C_F(m-\bar m)}{\tilde u}+
\dual{C_G(m-\bar m)}{\tilde p}\right\}\dd t
+\gamma(\bar m,m-\bar m)_\UU+\mathcal E_0(m-\bar m)\ge0
\label{eq:full-opt-condition}
\end{equation}
for all $m\in\UU_{\ad}$. This condition is necessary and sufficient. In particular, for a distributed fluid-source control $G=G_b+m_p$, that is, $C_F=0$, $C_G$ the canonical embedding, $m=m_p$, and $\UU=L^2H$, and for fixed initial data,
\begin{equation}
(\tilde p+\gamma\bar m_p,m_p-\bar m_p)_{L^2H}\ge0
\quad(m_p\in\UU_{\ad}).
\label{eq:source-control-VI}
\end{equation}
If $\UU_{\ad}=\{m_p:m_a\le m_p\le m_b\text{ a.e.}\}$ with measurable lower and upper bounds, then
\begin{equation}
\bar m_p(x,t)=\Pi_{[m_a(x,t),m_b(x,t)]}\!\left(-\frac{\tilde p(x,t)}\gamma\right)
\quad\text{for a.e. }(x,t).
\label{eq:source-proj}
\end{equation}
\end{theorem}
\begin{proof}
The derivative of \eqref{eq:cost} in direction $\eta$ is
\begin{equation*}
j'(\bar m)\eta=\int_0^T\!\left\{(R_u,w)_{L^2(\Omega)^d}+(R_p,r)_H\right\}\dd t+\gamma(\bar m,\eta)_\UU,
\end{equation*}
where $(w,r)$ is the state variation generated by $\eta$. Insert Proposition~\ref{prop:duality} with $\eta=m-\bar m$ to obtain \eqref{eq:full-opt-condition}. Theorem~\ref{thm:ocp-unique} makes this variational inequality necessary and sufficient. For $C_F=0$, $C_G$ the canonical embedding of $L^2H$ into $L^2Q^*$, and fixed independent initial data, the endpoint term is zero and \eqref{eq:source-control-VI} follows. The pointwise characterization is the standard orthogonal projection onto an $L^2$ box-constraint set.
\end{proof}
\begin{remark}[Boundary and mechanical controls]
For a normal Darcy flux with the outward-flux sign convention \eqref{eq:loads}, the fluid contribution to the gradient is $-\tilde p$ restricted to the flux boundary. For mechanical traction the time-integral contribution is the trace of $\tilde u$; if $\delta=0$ and the initial pressure is fixed, the endpoint term \eqref{eq:initial-correction} also contributes unless the initial traction is fixed. If mechanical controls are regularized in an $H^1$-in-time norm, their Riesz representatives must be taken in that control Hilbert space: the optimality condition is then generally \emph{not} an $L^2$ pointwise projection.
\end{remark}

\section{Singular limits for the state and adjoint equations}\label{sec:limits}
All constants $C$ in this section are generic in the sense of Section~\ref{subsec:mass}; in particular, they are uniform in
\begin{equation}
0\le c_0\le\bar c_0,\qquad 0\le\delta\le\bar\delta.
\label{eq:parameter-rectangle}
\end{equation}
We distinguish a stronger energy-data setting from a weaker one because different convergence norms require different initial regularity.

\subsection{Uniform estimates for the elastic family}
\begin{assumption}[Compatible data for the singular limits]\label{ass:regular-data}
The mechanical load and fluid source satisfy
\begin{equation}
F\in H^1V^*,\qquad G\in L^2H,\qquad p_0\in Q.
\label{eq:regular-data}
\end{equation}
For every $c_0,\delta$ we select the same initial pressure $p_0$ wherever a pressure initial datum is meaningful; at $\delta>0$, where the initial displacement is an independent datum, we choose it to be the elastic initial displacement \eqref{eq:elastic-compatibility}, which at $\delta=0$ is forced by the equations:
\begin{equation}
u_0=A_{\ee}^{-1}(F(0)+B^*p_0).
\label{eq:prepared}
\end{equation}
For $c_0=0$, $\delta>0$ this does not impose an additional independent pressure initial condition: it is only a preparation of the displacement against the target elastic equilibrium.
\end{assumption}

\begin{lemma}[Uniform time regularity of the elastic family]\label{lem:uniform-reg}
Under Assumption~\ref{ass:regular-data}, the purely elastic states $(u_{c_0,0},p_{c_0,0})$ satisfy
\begin{equation}
\begin{split}
\sup_{0\le c_0\le\bar c_0}\Bigl[
\norm{\partial_t p_{c_0,0}}_{L^2H}
+\norm{\partial_t u_{c_0,0}}_{L^2V}
+\norm{p_{c_0,0}}_{L^\infty Q}\Bigr]
\le C\mathcal D,
\end{split}
\label{eq:uniform-reg}
\end{equation}
where
\begin{equation}
\mathcal D:=\norm{p_0}_Q+\norm G_{L^2H}
+\norm{F_t}_{L^2V^*}.
\label{eq:D}
\end{equation}
\end{lemma}
\begin{proof}
Write $G_F:=G-BA_{\ee}^{-1}F_t$ as in the proof of Theorem~\ref{thm:state-storage}. The reduced pressure equation is
\begin{equation}
(M_{c_0}p_{c_0,0})_t+\KK p_{c_0,0}=G_F\quad\text{in }Q^*,
\label{eq:uniform-regular-equation}
\end{equation}
where $G_F\in L^2H$ and $\norm{G_F}_{L^2H}\le C(\norm G_{L^2H}+\norm{F_t}_{L^2V^*})$. Let $p_n\in H^1(0,T;Q_n)$ be the Galerkin approximation in the spaces $Q_n$ of the proof of Theorem~\ref{thm:state-full}:
\begin{equation*}
(M_{c_0}\partial_tp_n,q)_H+a_p(t;p_n,q)=(G_F,q)_H\quad\text{for all }q\in Q_n,\qquad p_n(0)=\Pi_np_0.
\end{equation*}
Since $\norm\cdot_Q$ is the form norm of $\KK_H(0)$, $\Pi_np_0\to p_0$ in $Q$ and $\norm{\Pi_np_0}_Q\le\norm{p_0}_Q$. Taking $q=\partial_tp_n$ gives
\begin{equation}
(M_{c_0} \partial_t p_n,\partial_t p_n)_H
+\frac12\frac{\dd}{\dd t}a_p\bigl(t;p_n(t),p_n(t)\bigr)
-\frac12a_p'(p_n,p_n)
=(G_F,\partial_t p_n)_H.
\label{eq:uniform-reg-identity}
\end{equation}
By Lemma~\ref{lem:mass-coercivity}, $(M_{c_0} x,x)_H\ge c_*\norm x_H^2$ independently of $c_0$. Apply Young's inequality, bound $\abs{a_p'(p_n,p_n)}\le\kappa_1\underline C_p^{-1}a_p(p_n,p_n)$, integrate over $(0,t)$, use $a_p(0;p_n(0),p_n(0))=\norm{\Pi_np_0}_Q^2\le\norm{p_0}_Q^2$, and apply Gronwall's inequality:
\begin{equation*}
\frac{c_*}2\norm{\partial_t p_n}_{L^2H}^2
+\frac12\sup_{0\le t\le T}a_p(p_n(t),p_n(t))
\le \exp(\kappa_1T/\underline C_p)\Bigl(\norm{p_0}_Q^2+\frac1{c_*}\norm{G_F}_{L^2H}^2\Bigr).
\end{equation*}
These bounds are uniform in $n$ and in $c_0$. Hence a subsequence of $(p_n)$ converges weakly in $H^1H$ and weakly-$*$ in $L^\infty Q$; its limit solves \eqref{eq:pressure-reduced} with $p(0)=p_0$ and has the regularity \eqref{eq:pressure-state-space}, so it equals $p_{c_0,0}$ by the uniqueness in Theorem~\ref{thm:state-storage}. By weak and weak-$*$ lower semicontinuity, the estimate holds for $p_{c_0,0}$, which therefore belongs to $H^1H\cap L^\infty Q$. Differentiating the \emph{representation} of the displacement, not the original mechanical equation as an independent PDE, gives
\begin{equation*}
\partial_t u_{c_0,0}=A_{\ee}^{-1}(F_t+B^*\partial_t p_{c_0,0})\in L^2V
\end{equation*}
with the stated uniform bound. Here the time independence of $A_{\ee}$ is used; the time dependence of $\mathbf K$ enters only through the factor $\exp(\kappa_1T/\underline C_p)$, which is independent of $c_0$.
\end{proof}

\subsection{Vanishing viscosity in the energy norm}
Fix $c_0\in[0,\bar c_0]$. For $\delta>0$ denote
\begin{equation}
e_u=u_{c_0,\delta}-u_{c_0,0},\qquad e_p=p_{c_0,\delta}-p_{c_0,0}.
\label{eq:err-delta}
\end{equation}
Subtract the equations and split $\partial_t u_{c_0,\delta}=\partial_t e_{u}+\partial_t u_{c_0,0}$ to find
\begin{align}
\delta a_{\vv}(\partial_t e_{u},v)+a_{\ee}(e_u,v)-(e_p,Bv)_H
&=-\delta a_{\vv}(\partial_t u_{c_0,0},v),\qquad v\in V,
\label{eq:delta-err-1}\\
(c_0e_p+Be_u)_t+\KK e_p&=0\quad\text{in }Q^*.
\label{eq:delta-err-2}
\end{align}
Preparation \eqref{eq:prepared} implies $e_u(0)=0$, and $c_0e_p(0)=0$ when $c_0>0$.

\begin{theorem}[Vanishing-viscosity estimate in the energy norm]\label{thm:visc-energy-rate}
Under Assumption~\ref{ass:regular-data}, uniformly in \eqref{eq:parameter-rectangle} with $\delta>0$,
\begin{equation}
\norm{e_u}_{L^\infty V}
+\norm{e_p}_{L^2Q}
+\sqrt\delta\norm{\partial_t e_{u}}_{L^2V}
\le C\sqrt\delta\,\mathcal D.
\label{eq:visc-energy-rate}
\end{equation}
In particular,
\begin{equation}
\norm{\partial_t u_{c_0,\delta}}_{L^2V}\le C\mathcal D
\label{eq:uniform-visc-velocity}
\end{equation}
with a constant independent of $c_0,\delta$.
\end{theorem}
\begin{proof}
In the error system \eqref{eq:delta-err-1}--\eqref{eq:delta-err-2}, test by $v=\partial_t e_{u}$ and $q=e_p$, as in the proof of Lemma~\ref{lem:state-energy}; this is admissible since $e_u\in H^1V$ by Theorems~\ref{thm:state-full}, \ref{thm:state-visc} and Lemma~\ref{lem:uniform-reg}, $e_p\in L^2Q$, and $c_0e_p\in H^1Q^*$ for $c_0>0$. The pressure/displacement coupling cancels, including when $c_0=0$, in which case the terms with $c_0\norm{e_p}_H^2$ below are absent. The resulting identity is
\begin{equation}
\begin{split}
&\frac12\frac{\dd}{\dd t}
\left\{a_{\ee}(e_u,e_u)+c_0\norm{e_p}_H^2\right\}
+\delta a_{\vv}(\partial_t e_{u},\partial_t e_{u})+a_p(e_p,e_p)\\
&\hspace{30mm}=-\delta a_{\vv}(\partial_t u_{c_0,0},\partial_t e_{u}).
\end{split}
\label{eq:delta-energy-identity}
\end{equation}
Cauchy--Schwarz in the $a_{\vv}$ inner product and Young's inequality give
\begin{equation*}
\delta\abs{a_{\vv}(\partial_t u_{c_0,0},\partial_t e_{u})}
\le\frac\delta2a_{\vv}(\partial_t e_{u},\partial_t e_{u})
+\frac\delta2a_{\vv}(\partial_t u_{c_0,0},\partial_t u_{c_0,0}).
\end{equation*}
The initial energy difference is zero. Integrating from $0$ to $t$, and then taking the supremum over $t$ in the first term and $t=T$ in the others, yields
\begin{align}
&\sup_{t\in[0,T]}\left\{a_{\ee}(e_u(t),e_u(t))+c_0\norm{e_p(t)}_H^2\right\}
+\delta\int_0^T a_{\vv}(\partial_t e_{u},\partial_t e_{u})\dd t
\nonumber\\
&\hspace{22mm}+2\int_0^T a_p(e_p,e_p)\dd t
\le2\delta\int_0^T a_{\vv}(\partial_t u_{c_0,0},\partial_t u_{c_0,0})\dd t.
\label{eq:delta-energy-controlled}
\end{align}
The estimate follows from the coercivity of the three forms and Lemma~\ref{lem:uniform-reg}. In particular, divide the velocity term by $\delta>0$ to obtain
$\norm{\partial_t e_{u}}_{L^2V}\le C\norm{\partial_t u_{c_0,0}}_{L^2V}$, with a constant independent of $\delta$. Since $\partial_t u_{c_0,\delta}=\partial_t u_{c_0,0}+\partial_t e_{u}$, \eqref{eq:uniform-visc-velocity} follows.
\end{proof}

\begin{theorem}[Optimality of the exponent $1/2$]\label{thm:visc-sharp}
Let $\mathcal C_{\rm v}=\mathcal C_{\rm e}$, $c_0\in[0,\bar c_0]$, and $0<\delta\le\min\{\bar\delta,T\}$. Then there are data satisfying Assumption~\ref{ass:regular-data}, with $G=0$, $p_0=0$, and $\mathcal D>0$, for which $e_p=0$ and
\begin{equation}
\norm{e_u}_{L^\infty V}\ge\frac{1-\exp(-1)}{\overline C_{\ee}}\sqrt\delta\,\mathcal D,
\qquad
\sqrt\delta\norm{\partial_t e_{u}}_{L^2V}\ge\frac{(1-\exp(-2))^{1/2}}{\sqrt2\,\overline C_{\ee}}\sqrt\delta\,\mathcal D.
\label{eq:visc-sharp}
\end{equation}
Consequently, for these two terms the exponent $1/2$ in \eqref{eq:visc-energy-rate} cannot be replaced by any $\sigma>1/2$ with a constant independent of $\delta$ and of the data.
\end{theorem}
The data in Theorem~\ref{thm:visc-sharp} may depend on $\delta$. The proof is given in Appendix~\ref{app:mechanical-sharpness}.

\subsection{Vanishing viscosity in the tracking norm}
Replacing $L^\infty V$ by $L^2V$ for the displacement and $L^2Q$ by $L^2H$ for the pressure gives a first-order estimate by backward duality.

\begin{theorem}[Vanishing-viscosity estimate in the tracking norm]\label{thm:visc-dual-rate}
Under Assumption~\ref{ass:regular-data}, there holds
\begin{equation}
\norm{u_{c_0,\delta}-u_{c_0,0}}_{L^2V}
+\norm{p_{c_0,\delta}-p_{c_0,0}}_{L^2H}
\le C\delta\,\mathcal D,
\label{eq:visc-Odelta}
\end{equation}
uniformly for $0\le c_0\le\bar c_0$ and $0\le\delta\le\bar\delta$.
\end{theorem}
\begin{proof}
Fix arbitrary dual observation data
\begin{equation*}
R_u\in L^2V^*,\qquad R_p\in L^2H.
\end{equation*}
Let $(\tilde u,\tilde p)$ solve the \emph{zero-viscosity} backward adjoint problem at the same value of $c_0$:
\begin{equation}
A_{\ee}\tilde u-B^*\tilde p_t=R_u,\qquad
-B\tilde u-c_0\tilde p_t+\KK \tilde p=R_p,\qquad \tilde p(T)=0.
\label{eq:dual-for-visc}
\end{equation}
By Theorem~\ref{thm:adj-elastic},
\begin{equation}
\norm {\tilde u}_{L^2V}\le C\left(\norm{R_u}_{L^2V^*}+\norm{R_p}_{L^2H}\right),
\label{eq:m-dual-uniform}
\end{equation}
where the constant does not depend on $c_0$.

Use the error equations \eqref{eq:delta-err-1}--\eqref{eq:delta-err-2}. Put $\Xi=c_0e_p+Be_u$. By \eqref{eq:delta-err-2}, $\Xi_t+\KK e_p=0$ in $Q^*$, and $\Xi(0)=0$ by the preparation \eqref{eq:prepared}. Moreover $\Xi\in C([0,T];H)$ and $\Xi_t\in L^2Q^*$, while $\tilde p\in H^1H\cap L^2Q$ by Theorem~\ref{thm:adj-elastic}, so Lemma~\ref{lem:product} applies to the pair $(\Xi,\tilde p)$. Insert \eqref{eq:dual-for-visc} into the observation pairing and use symmetry of $a_{\ee}$:
\begin{align*}
&\int_0^T\!\left\{\dual{R_u}{e_u}+(R_p,e_p)_H\right\}\dd t\\
&=\int_0^T\!\left\{a_{\ee}(\tilde u,e_u)-(\tilde p_t,Be_u)_H
-(B\tilde u,e_p)_H-c_0(\tilde p_t,e_p)_H+a_p(\tilde p,e_p)\right\}\dd t\\
&=\int_0^T\!\left\{\dual{A_{\ee}e_u-B^*e_p}{\tilde u}
-(\tilde p_t,\Xi)_H+a_p(\tilde p,e_p)\right\}\dd t\\
&=-\delta\int_0^T a_{\vv}(\partial_t u_{c_0,\delta},\tilde u)\dd t
-\int_0^T\!\left\{(\tilde p_t,\Xi)_H+\dual{\Xi_t}{\tilde p}\right\}\dd t.
\end{align*}
The last integral is $-[ (\Xi,\tilde p)_H]_0^T$ by Lemma~\ref{lem:product}; it vanishes because $\Xi(0)=0$ and $\tilde p(T)=0$. Consequently the exact identity
\begin{equation}
\int_0^T\!\left\{\dual{R_u}{e_u}+(R_p,e_p)_H\right\}\dd t
=-\delta\int_0^T a_{\vv}(\partial_t u_{c_0,\delta},\tilde u)\dd t
\label{eq:delta-duality-exact}
\end{equation}
holds.

Now apply the continuity of $a_{\vv}$, \eqref{eq:uniform-visc-velocity}, and \eqref{eq:m-dual-uniform}:
\begin{align*}
\abs{\int_0^T\!\left\{\dual{R_u}{e_u}+(R_p,e_p)_H\right\}\dd t}
&\le C\delta\norm{\partial_t u_{c_0,\delta}}_{L^2V}\norm {\tilde u}_{L^2V}\\
&\le C\delta\,\mathcal D
\left(\norm{R_u}_{L^2V^*}+\norm{R_p}_{L^2H}\right).
\end{align*}
Taking the supremum over the unit ball of $\YY_V^*=L^2V^*\times L^2H$ gives \eqref{eq:visc-Odelta}.
\end{proof}

\begin{proposition}[Sharpness of the pressure rates for fluid sources]\label{prop:source-sharp}
There is a configuration on $\Omega=(0,1)^2$ satisfying the assumptions of Section~\ref{sec:setting}, with $c_0=0$ and constant anisotropic tensors $\mathcal C_{\rm e}=\mathcal C_{\rm v}$, for which the following holds. For each sufficiently small $\delta>0$, there is a distributed source $G_\delta$ of unit norm in $L^2H$. With this same source, $F=0$, $p_0=0$, and zero boundary flux, the viscous and elastic pressures satisfy
\begin{equation*}
\norm{p_{0,\delta}-p_{0,0}}_{L^2H}\ge c\,\delta,\qquad
\norm{p_{0,\delta}-p_{0,0}}_{L^2Q}\ge c\,\sqrt\delta,
\end{equation*}
where $c>0$ is independent of $\delta$. These data satisfy Assumption~\ref{ass:regular-data} with $\mathcal D=1$. Thus neither pressure exponent in Theorems~\ref{thm:visc-energy-rate} and \ref{thm:visc-dual-rate} can be improved uniformly over distributed sources.
\end{proposition}
The configuration and the proof of Proposition~\ref{prop:source-sharp} are given in Appendix~\ref{subsec:example}.

\subsection{Vanishing storage}\label{subsec:storage}
Set $\delta=0$, and denote by $(u_{c_0,0},p_{c_0,0})$ and $(u_{0,0},p_{0,0})$ the solutions at positive storage $c_0$ and zero storage, with the same data and the same initial pressure. We write
\begin{equation*}
E_u:=u_{c_0,0}-u_{0,0},\qquad E_p:=p_{c_0,0}-p_{0,0}.
\end{equation*}

\begin{theorem}[Vanishing-storage estimate in the energy norm]\label{thm:storage-energy}
Under Assumption~\ref{ass:regular-data},
\begin{equation}
\norm{p_{c_0,0}-p_{0,0}}_{L^\infty H}
+\norm{p_{c_0,0}-p_{0,0}}_{L^2Q}
+\norm{u_{c_0,0}-u_{0,0}}_{L^\infty V}
\le Cc_0\,\mathcal D.
\label{eq:storage-energy-rate}
\end{equation}
The constant is uniform for $0\le c_0\le\bar c_0$.
\end{theorem}
\begin{proof}
Under Assumption~\ref{ass:regular-data}, Lemma~\ref{lem:uniform-reg} gives $\partial_tp_{c_0,0},\partial_tp_{0,0}\in L^2H$, hence $E_p\in H^1H\cap L^\infty Q$. Subtract the reduced pressure equations \eqref{eq:pressure-reduced} at $c_0$ and at $0$, which now hold with $(M_{c_0}p)_t=M_{c_0}\partial_tp\in L^2H$, and use $M_{c_0}=M_0+c_0I$:
\begin{equation}
M_{c_0}\partial_tE_p+\KK E_p=-c_0\,\partial_t p_{0,0},\qquad E_p(0)=0.
\label{eq:storage-error}
\end{equation}
Testing \eqref{eq:storage-error} directly with $E_p$ and using the symmetry and time independence of $M_{c_0}$, we obtain
\begin{equation*}
\frac12\frac{\dd}{\dd t}(M_{c_0}E_p,E_p)_H+a_p(E_p,E_p)
=-c_0(\partial_t p_{0,0},E_p)_H.
\end{equation*}
By Poincar\'e's inequality, boundedness of the embedding $Q\hookrightarrow H$, and Young's inequality,
\begin{equation*}
c_0\abs{(\partial_t p_{0,0},E_p)_H}
\le\frac12a_p(E_p,E_p)+Cc_0^2\norm{\partial_t p_{0,0}}_H^2.
\end{equation*}
Because $E_p(0)=0$, integration gives
\begin{equation*}
\sup_t(M_{c_0}E_p(t),E_p(t))_H+\norm{E_p}_{L^2Q}^2
\le Cc_0^2\norm{\partial_t p_{0,0}}_{L^2H}^2.
\end{equation*}
The uniform coercivity \eqref{eq:uniform-mass} bounds the $L^\infty H$ norm. The displacement difference is exactly
\begin{equation}
E_u=A_{\ee}^{-1}B^*E_p,
\label{eq:storage-displacement-relation}
\end{equation}
so boundedness of $A_{\ee}^{-1}B^*:H\to V$ yields the stated $L^\infty V$ estimate.
\end{proof}

A duality estimate requires less temporal regularity of the state pressure and yields an $O(c_0)$ tracking rate for energy-level initial data $p_0\in H$.
\begin{theorem}[Vanishing-storage estimate in the tracking norm]\label{thm:storage-duality}
Suppose only
\begin{equation}
F\in H^1V^*,\qquad G\in L^2Q^*,\qquad p_0\in H,
\label{eq:rough-storage-data}
\end{equation}
and take the same initial pressure $p_0$ for $c_0>0$ and $c_0=0$. Then
\begin{equation}
\norm{u_{c_0,0}-u_{0,0}}_{L^2V}
+\norm{p_{c_0,0}-p_{0,0}}_{L^2H}
\le Cc_0\,\mathcal D_0,
\label{eq:storage-L2-rate}
\end{equation}
where
\begin{equation}
\mathcal D_0:=\norm{p_0}_H+\norm G_{L^2Q^*}+\norm{F_t}_{L^2V^*},
\label{eq:rough-D}
\end{equation}
and $C$ is independent of $c_0\in[0,\bar c_0]$.
\end{theorem}
\begin{proof}
Let $E_u$ and $E_p$ be the storage differences defined at the beginning of Section~\ref{subsec:storage}. By \eqref{eq:u-from-p}, $A_{\ee}E_u-B^*E_p=0$, and $E_u(0)=0$ because both states have the same initial pressure and mechanical load. Subtracting the fluid equations \eqref{eq:state-fluid} of the two states shows that
\begin{equation}
\Xi:=c_0p_{c_0,0}+BE_u
\quad\text{satisfies}\quad
\Xi_t+\KK E_p=0\ \text{in }L^2Q^*,\qquad \Xi(0)=c_0p_0.
\label{eq:storage-difference-mass}
\end{equation}
By Theorem~\ref{thm:state-storage}, $\Xi\in C([0,T];H)$. Fix arbitrary $R_u\in L^2V^*$ and $R_p\in L^2H$ and solve the \emph{zero-storage, zero-viscosity} adjoint
\begin{equation}
A_{\ee}\tilde u-B^*\tilde p_t=R_u,\qquad
-B\tilde u+\KK \tilde p=R_p,\qquad \tilde p(T)=0.
\label{eq:storage-dual}
\end{equation}
Theorem~\ref{thm:adj-elastic} implies $\tilde p\in H^1H\cap L^2Q$ and
\begin{equation}
\norm{\tilde p_t}_{L^2H}+\norm{\tilde p(0)}_H
\le C\left(\norm{R_u}_{L^2V^*}+\norm{R_p}_{L^2H}\right),
\label{eq:h-test-bound}
\end{equation}
since $\tilde p(T)=0$ and $\norm{\tilde p(0)}_H\le\sqrt T\norm{\tilde p_t}_{L^2H}$. Since $A_{\ee}E_u-B^*E_p=0$, \eqref{eq:storage-difference-mass}--\eqref{eq:storage-dual} give
\begin{align*}
\int_0^T\!\left\{\dual{R_u}{E_u}+(R_p,E_p)_H\right\}\dd t
&=\int_0^T\!\left\{-(\tilde p_t,BE_u)_H+a_p(\tilde p,E_p)\right\}\dd t\\
&=\int_0^T\!\left\{-(\tilde p_t,\Xi)_H-\dual{\Xi_t}{\tilde p}+c_0(\tilde p_t,p_{c_0,0})_H\right\}\dd t.
\end{align*}
Lemma~\ref{lem:product} applies to $\Xi$ and $\tilde p$ and gives $-\int_0^T\{(\tilde p_t,\Xi)_H+\dual{\Xi_t}{\tilde p}\}\dd t=(\Xi(0),\tilde p(0))_H=c_0(p_0,\tilde p(0))_H$. Hence
\begin{equation}
\int_0^T\!\left\{\dual{R_u}{E_u}+(R_p,E_p)_H\right\}\dd t
=c_0(p_0,\tilde p(0))_H+c_0\int_0^T(p_{c_0,0},\tilde p_t)_H\dd t.
\label{eq:storage-transposition-identity}
\end{equation}
Theorem~\ref{thm:state-storage} yields $\norm{p_{c_0,0}}_{L^2H}\le C\mathcal D_0$ uniformly in $c_0$. Therefore the right-hand side is bounded by
\begin{equation*}
Cc_0\mathcal D_0\left(\norm{R_u}_{L^2V^*}+\norm{R_p}_{L^2H}\right).
\end{equation*}
Duality in $L^2V\times L^2H$ proves \eqref{eq:storage-L2-rate}.
\end{proof}

\subsection{The joint singular limit}
\begin{corollary}[Joint singular-limit estimates]\label{cor:combined-rate}
Under Assumption~\ref{ass:regular-data}, the state family satisfies
\begin{equation}
\norm{u_{c_0,\delta}-u_{0,0}}_{L^2V}
+\norm{p_{c_0,\delta}-p_{0,0}}_{L^2H}
\le C(c_0+\delta)\mathcal D.
\label{eq:combined-Oc-delta}
\end{equation}
In the energy norm one has the distinct, generally weaker bound
\begin{equation}
\norm{u_{c_0,\delta}-u_{0,0}}_{L^\infty V}
+\norm{p_{c_0,\delta}-p_{0,0}}_{L^2Q}
\le C(c_0+\sqrt\delta)\mathcal D.
\label{eq:combined-energy}
\end{equation}
In particular, both estimates imply convergence in the observation space $\YY$ defined in \eqref{eq:tracking-spaces}.
\end{corollary}
\begin{proof}
Insert the intermediate state $(u_{c_0,0},p_{c_0,0})$ and write
\begin{equation*}
(u_{c_0,\delta},p_{c_0,\delta})-(u_{0,0},p_{0,0})
=\bigl[(u_{c_0,\delta},p_{c_0,\delta})-(u_{c_0,0},p_{c_0,0})\bigr]
 +\bigl[(u_{c_0,0},p_{c_0,0})-(u_{0,0},p_{0,0})\bigr].
\end{equation*}
In $\YY_V$ apply Theorems~\ref{thm:visc-dual-rate} and \ref{thm:storage-duality}; the latter applies because Assumption~\ref{ass:regular-data} implies \eqref{eq:rough-storage-data} with $\mathcal D_0\le C\mathcal D$, since $\norm{p_0}_H\le C\norm{p_0}_Q$ and $\norm G_{L^2Q^*}\le C\norm G_{L^2H}$. For the stronger norm use Theorems~\ref{thm:visc-energy-rate} and \ref{thm:storage-energy}. The constants are uniform in $c_0$ because all uses of the pressure-to-dilation mass are controlled by \eqref{eq:uniform-mass}.
\end{proof}
\begin{remark}[Uniformity of the estimates]
Under Assumption~\ref{ass:regular-data}, the constants in \eqref{eq:combined-Oc-delta}--\eqref{eq:combined-energy} are independent of $c_0$ and $\delta$. These estimates do not give uniform time regularity of the adjoint or cover arbitrary mechanical loads $F\in L^2V^*$ as $\delta\downarrow0$. Boundary fluxes generally give only $G\in L^2Q^*$; their convergence is treated in Section~\ref{sec:boundary-controls}.
\end{remark}

\subsection{Uniform bounds and convergence of the adjoint variables}\label{subsec:adjoint-limit}
For $G\in L^2Q^*$ let $\mathcal T_{c_0,\delta}G:=(w,r)$ denote the state with $F=0$, fluid load $G$, and homogeneous initial data, given by Theorems~\ref{thm:state-full}, \ref{thm:state-visc}, and \ref{thm:state-storage}. For $\delta=0$ the homogeneous initial data mean $r(0)=0$, and then $w(0)=0$ by \eqref{eq:elastic-compatibility}.

\begin{lemma}[Uniform stability for fluid loads]\label{lem:Qstar-stability}
Uniformly for $(c_0,\delta)$ in \eqref{eq:parameter-rectangle},
\begin{equation}
\norm w_{L^2V}+\norm r_{L^2Q}\le C\norm G_{L^2Q^*},\qquad (w,r)=\mathcal T_{c_0,\delta}G.
\label{eq:Qstar-stability}
\end{equation}
\end{lemma}
\begin{proof}
For $\delta>0$, Lemma~\ref{lem:state-energy} with $F=0$ and Young's inequality give
\begin{equation*}
\frac{\dd}{\dd t}\bigl\{a_{\ee}(w,w)+c_0\norm r_H^2\bigr\}+a_p(r,r)
\le\underline C_p^{-1}\norm G_{Q^*}^2,
\end{equation*}
and integration with zero initial data bounds $\norm w_{L^\infty V}+\norm r_{L^2Q}$ by $C\norm G_{L^2Q^*}$; the $\delta$-dependence of the constant in \eqref{eq:visc-est} arises only from $F$. For $\delta=0$, \eqref{eq:pressure-state-est} with $F=0$ and $r(0)=0$ bounds $\norm r_{L^2Q}$, and $w=A_{\ee}^{-1}B^*r$ bounds $\norm w_{L^2V}$.
\end{proof}

\begin{lemma}[Multiplicative trace inequality]\label{lem:mult-trace}
Let $\mathrm{Tr}_p:Q\to L^2(\Gamma_{p,N})$ denote the trace operator. There is $C_{\rm tr}$ such that
\begin{equation}
\norm{\mathrm{Tr}_pq}_{L^2(\Gamma_{p,N})}^2\le C_{\rm tr}\norm q_H\norm q_Q,
\qquad q\in Q.
\label{eq:mult-trace}
\end{equation}
\end{lemma}
\begin{proof}
The estimate follows from \cite[Theorem~1.5.1.10]{Grisvard1985} with $p=2$ and $\varepsilon=\norm q_{L^2(\Omega)}/\norm q_{H^1(\Omega)}$, by restriction to $\Gamma_{p,N}$ and the equivalence of $\norm\cdot_Q$ and $\norm\cdot_{H^1(\Omega)}$ on $Q$.
\end{proof}

\begin{theorem}[Adjoint variables in the joint limit]\label{thm:adjoint-limit}
Let $R_u\in L^2V^*$ and $R_p\in L^2H$, put $\mathcal R_*:=\norm{R_u}_{L^2V^*}+\norm{R_p}_{L^2H}$, and for $(c_0,\delta)$ in \eqref{eq:parameter-rectangle} let $(\tilde u_{c_0,\delta},\tilde p_{c_0,\delta})$ be the adjoint pair of Theorem~\ref{thm:adj-full}, \ref{thm:adj-elastic}, or \ref{thm:adj-visc} with these data. Then, uniformly in $(c_0,\delta)$,
\begin{align}
\norm{\tilde u_{c_0,\delta}}_{L^2V}+\norm{\tilde p_{c_0,\delta}}_{L^2Q}&\le C\mathcal R_*,
\label{eq:adj-uniform-Q}\\
\norm{\tilde p_{c_0,\delta}-\tilde p_{c_0,0}}_{L^2H}\le C\delta\mathcal R_*,\qquad
\norm{\tilde p_{c_0,\delta}-\tilde p_{c_0,0}}_{L^2Q}&\le C\sqrt\delta\,\mathcal R_*,
\label{eq:adj-rates-visc}\\
\norm{\tilde u_{c_0,0}-\tilde u_{0,0}}_{L^2V}
+\norm{\partial_t(\tilde p_{c_0,0}-\tilde p_{0,0})}_{L^2H}
+\norm{\tilde p_{c_0,0}-\tilde p_{0,0}}_{L^\infty Q}&\le Cc_0\mathcal R_*,
\label{eq:adj-rates}
\end{align}
and consequently
\begin{equation}
\begin{gathered}
\norm{\tilde p_{c_0,\delta}-\tilde p_{0,0}}_{L^2H}\le C(c_0+\delta)\mathcal R_*,\qquad
\norm{\tilde p_{c_0,\delta}-\tilde p_{0,0}}_{L^2Q}\le C(c_0+\sqrt\delta)\mathcal R_*,\\
\norm{\mathrm{Tr}_p(\tilde p_{c_0,\delta}-\tilde p_{0,0})}_{L^2L^2(\Gamma_{p,N})}\le C\bigl(c_0+\delta^{3/4}\bigr)\mathcal R_*.
\end{gathered}
\label{eq:adj-joint}
\end{equation}
Moreover, $\tilde u_{c_0,\delta}\to\tilde u_{0,0}$ in $L^2V$ as $(c_0,\delta)\to(0,0)$.
\end{theorem}
\begin{proof}
\emph{Step 1: vanishing viscosity in $L^2H$ by duality.} For $G\in L^2Q^*$ and $(w,r)=\mathcal T_{c_0,\delta}G$, Proposition~\ref{prop:duality} with $(F_w,G_w)=(0,G)$ and zero initial variations gives $\mathcal E_0=0$ and
\begin{equation}
\int_0^T\dual G{\tilde p_{c_0,\delta}}\dd t
=\int_0^T\!\left\{\dual{R_u}w+(R_p,r)_H\right\}\dd t.
\label{eq:adj-duality-G}
\end{equation}
Subtract \eqref{eq:adj-duality-G} at $(c_0,\delta)$ and at $(c_0,0)$. For $G\in L^2H$ the data $(F,G,p_0)=(0,G,0)$ satisfy Assumption~\ref{ass:regular-data} with $\mathcal D=\norm G_{L^2H}$ and prepared displacement $u_0=0$, so Theorem~\ref{thm:visc-dual-rate} gives
\begin{equation*}
\abs{\int_0^T(G,\tilde p_{c_0,\delta}-\tilde p_{c_0,0})_H\dd t}
\le C\mathcal R_*\norm{\mathcal T_{c_0,\delta}G-\mathcal T_{c_0,0}G}_{\YY_V}
\le C\delta\mathcal R_*\norm G_{L^2H}.
\end{equation*}
Taking the supremum over $G$ in the unit ball of $L^2H$ gives the first bound in \eqref{eq:adj-rates-visc}.

\emph{Step 2: an integrated error equation and uniform stability.} Let $\delta>0$ and put $e_{\tilde u}:=\tilde u_{c_0,\delta}-\tilde u_{c_0,0}$, $e_{\tilde p}:=\tilde p_{c_0,\delta}-\tilde p_{c_0,0}$, and
\begin{equation*}
h(t):=\int_t^Te_{\tilde u}(s)\dd s,\qquad\text{so that}\qquad h\in H^1V,\quad h_t=-e_{\tilde u},\quad h(T)=0.
\end{equation*}
Subtracting the adjoint equations \eqref{eq:adj-mech}--\eqref{eq:adj-fluid} at $(c_0,\delta)$ and $(c_0,0)$ gives
\begin{equation*}
-\Lambda_t+A_{\ee}e_{\tilde u}=0,\qquad \Lambda:=\delta A_{\vv}\tilde u_{c_0,\delta}+B^*e_{\tilde p},
\qquad\quad
-(c_0e_{\tilde p})_t-Be_{\tilde u}+\KK e_{\tilde p}=0.
\end{equation*}
Here $\Lambda=z_{c_0,\delta}-B^*\tilde p_{c_0,0}\in H^1V^*$, where $z_{c_0,\delta}$ is the quantity $z$ of Section~\ref{subsec:adj-form} for $(\tilde u_{c_0,\delta},\tilde p_{c_0,\delta})$, and $\Lambda(T)=0$ by \eqref{eq:terminal-general} and $\tilde p_{c_0,0}(T)=0$. Integrating $-\Lambda_t+A_{\ee}e_{\tilde u}=0$ over $(t,T)$ gives $\Lambda=-A_{\ee}h$. Since $\tilde u_{c_0,\delta}=\tilde u_{c_0,0}-h_t$, the pair $(h,e_{\tilde p})$ satisfies
\begin{subequations}\label{eq:adj-integrated-error}
\begin{align}
-\delta A_{\vv}h_t+A_{\ee}h+B^*e_{\tilde p}&=-\delta A_{\vv}\tilde u_{c_0,0}\quad\text{in }L^2V^*,
\label{eq:adj-integrated-error-mech}\\
-(c_0e_{\tilde p})_t+Bh_t+\KK e_{\tilde p}&=0\quad\text{in }L^2Q^*.
\label{eq:adj-integrated-error-fluid}
\end{align}
\end{subequations}
Test \eqref{eq:adj-integrated-error-mech} by $-h_t\in L^2V$ and \eqref{eq:adj-integrated-error-fluid} by $e_{\tilde p}\in L^2Q$. For $c_0>0$ both pressure adjoints lie in $H^1H$, and for $c_0=0$ the time-derivative term is absent. The coupling terms $-(e_{\tilde p},Bh_t)_H$ and $(Bh_t,e_{\tilde p})_H$ cancel. Since $h(T)=0$ and, for $c_0>0$, $e_{\tilde p}(T)=0$, integration over $(t,T)$ gives
\begin{equation}
\begin{aligned}
&\frac12a_{\ee}(h(t),h(t))+\frac{c_0}2\norm{e_{\tilde p}(t)}_H^2
+\int_t^T\!\Bigl\{\delta a_{\vv}(h_t,h_t)+a_p(s;e_{\tilde p},e_{\tilde p})\Bigr\}\dd s\\
&\qquad=\delta\int_t^Ta_{\vv}(\tilde u_{c_0,0},h_t)\dd s.
\end{aligned}
\label{eq:adj-integrated-energy}
\end{equation}
By Young's inequality on the right-hand side, \eqref{eq:ellipticity}, and the $c_0$-uniform bound $\norm{\tilde u_{c_0,0}}_{L^2V}\le C\mathcal R_*$ of Theorem~\ref{thm:adj-elastic},
\begin{equation}
\norm h_{L^\infty V}+\sqrt\delta\norm{h_t}_{L^2V}+\norm{e_{\tilde p}}_{L^2Q}\le C\sqrt\delta\,\mathcal R_*.
\label{eq:adj-h-bound}
\end{equation}
This is the second bound in \eqref{eq:adj-rates-visc}. Together with $\norm{\tilde u_{c_0,0}}_{L^2V}+\norm{\tilde p_{c_0,0}}_{L^\infty Q}\le C\mathcal R_*$ from Theorem~\ref{thm:adj-elastic} and $\norm{e_{\tilde u}}_{L^2V}=\norm{h_t}_{L^2V}\le C\mathcal R_*$, it gives \eqref{eq:adj-uniform-Q}.

\emph{Step 3: vanishing storage.} Put $\hat e:=\tilde p_{c_0,0}-\tilde p_{0,0}$. Subtracting \eqref{eq:adj-reduced-pressure} at $c_0$ and at $0$, and using $M_{c_0}=M_0+c_0I$,
\begin{equation*}
-M_{c_0}\hat e_t+\KK\hat e=c_0\,\partial_t\tilde p_{0,0},\qquad \hat e(T)=0.
\end{equation*}
This is \eqref{eq:adj-reduced-pressure} with data $(R_u,R_p)=(0,c_0\partial_t\tilde p_{0,0})$, and $\tilde u_{c_0,0}-\tilde u_{0,0}=A_{\ee}^{-1}B^*\hat e_t$ by \eqref{eq:adj-elim-m}. Theorem~\ref{thm:adj-elastic}, applied once with these data and once at $c_0=0$ with the original data, gives \eqref{eq:adj-rates}:
\begin{equation*}
\norm{\hat e_t}_{L^2H}+\norm{\hat e}_{L^\infty Q}+\norm{\tilde u_{c_0,0}-\tilde u_{0,0}}_{L^2V}
\le Cc_0\norm{\partial_t\tilde p_{0,0}}_{L^2H}\le Cc_0\mathcal R_*.
\end{equation*}

\emph{Step 4: joint limit and traces.} The first two bounds in \eqref{eq:adj-joint} follow from \eqref{eq:adj-rates-visc}--\eqref{eq:adj-rates} and the triangle inequality. By \eqref{eq:mult-trace}, the Cauchy--Schwarz inequality in time, and \eqref{eq:adj-rates-visc},
\begin{equation*}
\norm{\mathrm{Tr}_pe_{\tilde p}}_{L^2L^2(\Gamma_{p,N})}^2
\le C_{\rm tr}\norm{e_{\tilde p}}_{L^2H}\norm{e_{\tilde p}}_{L^2Q}
\le C\delta^{3/2}\mathcal R_*^2,
\end{equation*}
while $\norm{\mathrm{Tr}_p\hat e}_{L^2L^2(\Gamma_{p,N})}\le\norm{\mathrm{Tr}_p}\norm{\hat e}_{L^2Q}\le Cc_0\mathcal R_*$. The triangle inequality gives the trace bound in \eqref{eq:adj-joint}.

\emph{Step 5: convergence of the mechanical adjoint.} Fix the residuals $(R_u,R_p)$ and let $(c_0^k,\delta_k)\to(0,0)$. Denote by $h^k$ the function $h$ of Step~2 for the pair $(c_0^k,\delta_k)$; indices with $\delta_k=0$ need no argument, since then $h^k=0$. By \eqref{eq:adj-h-bound}, $h^k\to0$ in $L^\infty V$ and $h^k_t$ is bounded in $L^2V$. Every weak limit point of $(h^k_t)$ in $L^2V$ is the distributional time derivative of $\lim_kh^k=0$, hence $h^k_t\rightharpoonup0$ weakly in $L^2V$. Dividing \eqref{eq:adj-integrated-energy} at $t=0$ by $\delta_k$ and dropping nonnegative terms,
\begin{equation*}
\int_0^Ta_{\vv}(h^k_t,h^k_t)\dd t
\le\int_0^Ta_{\vv}(\tilde u_{0,0},h^k_t)\dd t
+\int_0^Ta_{\vv}(\tilde u_{c_0^k,0}-\tilde u_{0,0},h^k_t)\dd t.
\end{equation*}
The first term tends to zero by weak convergence, and the second is bounded by $Cc_0^k\mathcal R_*^2$ by \eqref{eq:adj-rates} and \eqref{eq:adj-h-bound}. Hence $h^k_t\to0$ in $L^2V$, and $\tilde u_{c_0^k,\delta_k}-\tilde u_{0,0}=-h^k_t+(\tilde u_{c_0^k,0}-\tilde u_{0,0})\to0$ in $L^2V$. Without further time regularity of $\tilde u_{0,0}$ this convergence comes without a rate, since the energy estimate \eqref{eq:adj-h-bound} controls $\tilde u_{c_0,\delta}-\tilde u_{c_0,0}$ in $L^2V$ only by a uniform bound; see Corollary~\ref{cor:mech-adjoint-rate}.
\end{proof}

\begin{corollary}[Rate for the mechanical adjoint]\label{cor:mech-adjoint-rate}
If, in addition, $\tilde u_{0,0}\in H^1V$, then, uniformly in \eqref{eq:parameter-rectangle},
\begin{equation}
\norm{\tilde u_{c_0,\delta}-\tilde u_{0,0}}_{L^2V}\le C\bigl(\sqrt\delta\,\norm{\tilde u_{0,0}}_{H^1V}+c_0\mathcal R_*\bigr).
\label{eq:mech-adjoint-rate}
\end{equation}
At $c_0=0$ the exponent $1/2$ in \eqref{eq:mech-adjoint-rate} is sharp even when the residuals are fixed and $\tilde u_{0,0}$ is constant in time; see Appendix~\ref{app:mechanical-sharpness}.
\end{corollary}
\begin{proof}
For $\delta=0$, \eqref{eq:mech-adjoint-rate} is contained in \eqref{eq:adj-rates}. Let $\delta>0$, let $h$ be as in Step~2 of the proof of Theorem~\ref{thm:adjoint-limit}, put $v:=\tilde u_{0,0}$, $X:=\norm h_{L^\infty V}$, and $Y:=\norm{h_t}_{L^2V}=\norm{\tilde u_{c_0,\delta}-\tilde u_{c_0,0}}_{L^2V}$, and split the right-hand side of \eqref{eq:adj-integrated-energy} with $\tilde u_{c_0,0}=v+(\tilde u_{c_0,0}-v)$. Since $a_{\vv}$ is independent of time and $h(T)=0$,
\begin{equation*}
\int_t^Ta_{\vv}(v,h_t)\dd s=-a_{\vv}\bigl(v(t),h(t)\bigr)-\int_t^Ta_{\vv}(v_t,h)\dd s,
\end{equation*}
whose modulus is at most $C\norm v_{H^1V}X$, while $\abs{\int_t^Ta_{\vv}(\tilde u_{c_0,0}-v,h_t)\dd s}\le Cc_0\mathcal R_*Y$ by \eqref{eq:adj-rates}. Taking the supremum over $t$ in the term $a_{\ee}(h(t),h(t))$ of \eqref{eq:adj-integrated-energy} and $t=0$ in the term $\delta\int_0^Ta_{\vv}(h_t,h_t)\dd s$ gives
\begin{equation*}
X^2+\delta Y^2\le C\delta\bigl(\norm v_{H^1V}X+c_0\mathcal R_*Y\bigr).
\end{equation*}
Young's inequality yields $X^2+\delta Y^2\le C(\delta^2\norm v_{H^1V}^2+\delta c_0^2\mathcal R_*^2)$, hence $Y\le C(\sqrt\delta\norm v_{H^1V}+c_0\mathcal R_*)$, and \eqref{eq:mech-adjoint-rate} follows from \eqref{eq:adj-rates} and the triangle inequality.
\end{proof}

\section{Convergence of constrained optimal controls}\label{sec:optimal-limits}
The parameter-dependent optimal control problem is controlled most conveniently through the operator perturbation of its reduced objective. For source and traction controls this route uses only the state estimates of Section~\ref{sec:limits}; the adjoint estimates of Section~\ref{subsec:adjoint-limit} enter only for flux controls.

\subsection{A perturbation result for strongly convex problems}
Let $\theta$ be an abstract parameter and suppose
\begin{equation}
j_\theta(m)=\frac12\norm{\mathcal S_\theta m+y_\theta}_{\YY}^2
+\frac\gamma2\norm m_\UU^2,
\qquad m\in\UU_{\ad},
\label{eq:abstract-perturbed-cost}
\end{equation}
where $\mathcal S_\theta\in\calL(\UU,\YY)$, $y_\theta\in\YY$ and $\gamma>0$ do not depend on $m$. Write $\Phi_\theta(m)=\mathcal S_\theta^*(\mathcal S_\theta m+y_\theta)$. The following estimate is the argument used in the error analysis of variationally discretized control problems \cite{Hinze2005}; we include the short proof for completeness.

\begin{lemma}[Perturbation of strongly convex minimization problems]\label{lem:perturbation}
Let $m_\theta,m_0$ be the unique minimizers of \eqref{eq:abstract-perturbed-cost} at $\theta$ and $0$ over the same nonempty closed convex set. Then
\begin{equation}
\norm{m_\theta-m_0}_\UU
\le\frac1\gamma\norm{\Phi_\theta(m_0)-\Phi_0(m_0)}_\UU.
\label{eq:convex-perturb-1}
\end{equation}
If $\sup_\theta\norm{\mathcal S_\theta}_{\calL(\UU,\YY)}\le C_{\mathcal S}$ and
\begin{equation}
\norm{\mathcal S_\theta-\mathcal S_0}_{\calL(\UU,\YY)}
+\norm{y_\theta-y_0}_\YY\le\omega_\theta,
\label{eq:convex-perturb-assumption}
\end{equation}
then
\begin{equation}
\norm{m_\theta-m_0}_\UU
\le\frac{C}{\gamma}\omega_\theta
\bigl(1+\norm{m_0}_\UU+\norm{y_0}_\YY\bigr),
\label{eq:convex-perturb-2}
\end{equation}
where $C$ depends only on $C_{\mathcal S}$ and $\norm{\mathcal S_0}$.
\end{lemma}
\begin{proof}
The two variational inequalities are
\begin{align*}
(\Phi_\theta(m_\theta)+\gamma m_\theta,m_0-m_\theta)_\UU&\ge0,\\
(\Phi_0(m_0)+\gamma m_0,m_\theta-m_0)_\UU&\ge0.
\end{align*}
Adding them and using $\Phi_\theta(m_\theta)-\Phi_\theta(m_0)=\mathcal S_\theta^*\mathcal S_\theta(m_\theta-m_0)$ gives
\begin{equation}
\gamma\norm{m_\theta-m_0}_\UU^2+\norm{\mathcal S_\theta(m_\theta-m_0)}_\YY^2
\le\bigl(\Phi_0(m_0)-\Phi_\theta(m_0),m_\theta-m_0\bigr)_\UU.
\label{eq:monotone-step}
\end{equation}
Dropping the second term on the left and applying the Cauchy--Schwarz inequality proves \eqref{eq:convex-perturb-1}. For the second estimate use the exact decomposition
\begin{align*}
\Phi_\theta(m_0)-\Phi_0(m_0)
&=(\mathcal S_\theta^*-\mathcal S_0^*)(\mathcal S_0m_0+y_0)\\
&\quad+\mathcal S_\theta^*(\mathcal S_\theta-\mathcal S_0)m_0
+\mathcal S_\theta^*(y_\theta-y_0).
\end{align*}
Its norm is bounded by a constant times
$\omega_\theta(1+\norm{m_0}_\UU+\norm{y_0}_\YY)$, using $\norm{\mathcal S_\theta^*-\mathcal S_0^*}=\norm{\mathcal S_\theta-\mathcal S_0}$. This proves \eqref{eq:convex-perturb-2}.
\end{proof}

\subsection{Distributed source controls}\label{subsec:source}
Let
\begin{equation}
\UU=L^2H,\quad \UU_{\ad}\subset\UU\text{ nonempty, closed, convex},
\quad F=F_b,\quad G=G_b+m_p.
\label{eq:source-family}
\end{equation}
This is \eqref{eq:affine-data} with $C_F=0$, $C_G$ the canonical embedding of $L^2H$ into $L^2Q^*$, and $m=m_p$. The initial data $u_0,p_0$ are independent of $m_p$. Assume
\begin{equation}
F_b\in H^1V^*,\quad G_b\in L^2H,\quad p_0\in Q,
\quad\text{and $u_0$ is given by \eqref{eq:prepared} with $F=F_b$};
\label{eq:ocp-prepared}
\end{equation}
that is, the background data $(F_b,G_b,p_0)$ satisfy Assumption~\ref{ass:regular-data}. These conditions ensure that Assumption~\ref{ass:regular-data} applies uniformly for each source control, with the same mechanical initial datum. In this setting the control-to-state map of \eqref{eq:control-map} is $\mathcal S_{c_0,\delta}=\mathcal T_{c_0,\delta}|_{L^2H}$, the state generated by the source $G=m_p$ with $F=0$ and homogeneous initial data, and $y^b_{c_0,\delta}$ is the state with loads $(F_b,G_b)$ and the initial data \eqref{eq:ocp-prepared}. We also regard $\mathcal S_{c_0,\delta}$ as a map into $\YY_V$; the reduced objective uses $\YY$, and $\YY_V$ is used only for estimates. We write
\begin{equation*}
\mathcal D_b:=\norm{p_0}_Q+\norm{G_b}_{L^2H}+\norm{\partial_tF_b}_{L^2V^*}
\end{equation*}
for the quantity \eqref{eq:D} of the background data.

\begin{proposition}[Convergence of the control-to-state operators]\label{prop:operator-rate}
Under \eqref{eq:ocp-prepared},
\begin{equation}
\sup_{0\le c_0\le\bar c_0,\,0\le\delta\le\bar\delta}
\norm{\mathcal S_{c_0,\delta}}_{\calL(\UU,\YY_V)}\le C,
\label{eq:uniform-G}
\end{equation}
and
\begin{equation}
\norm{\mathcal S_{c_0,\delta}-\mathcal S_{0,0}}_{\calL(\UU,\YY_V)}
+\norm{y^b_{c_0,\delta}-y^b_{0,0}}_{\YY_V}
\le C(c_0+\delta)(1+\mathcal D_b).
\label{eq:operator-rate}
\end{equation}
The same bounds hold after replacing $\YY_V$ by $\YY$.
\end{proposition}
\begin{proof}
For each $m_p\in\UU$, the homogeneous-state data satisfy Assumption~\ref{ass:regular-data} with $F=0$, $p_0=0$, and $\mathcal D=\norm{m_p}_{L^2H}$. Corollary~\ref{cor:combined-rate} therefore gives
\begin{equation*}
\norm{(\mathcal S_{c_0,\delta}-\mathcal S_{0,0})m_p}_{\YY_V}
\le C(c_0+\delta)\norm{m_p}_\UU.
\end{equation*}
Taking the supremum over $\norm{m_p}_\UU=1$ proves the first term of \eqref{eq:operator-rate}. The background states satisfy $\norm{y^b_{c_0,\delta}-y^b_{0,0}}_{\YY_V}\le C(c_0+\delta)\mathcal D_b$ by Corollary~\ref{cor:combined-rate}, because $F_b,G_b,p_0$ satisfy \eqref{eq:ocp-prepared}. Finally, boundedness of $\mathcal S_{0,0}$ follows from Theorem~\ref{thm:state-storage} and the energy estimate, while
\begin{equation*}
\norm{\mathcal S_{c_0,\delta}}\le\norm{\mathcal S_{0,0}}+\norm{\mathcal S_{c_0,\delta}-\mathcal S_{0,0}}
\end{equation*}
proves \eqref{eq:uniform-G}. The embedding $\YY_V\hookrightarrow\YY$ is continuous, giving the asserted tracking-space version.
\end{proof}

\begin{theorem}[Convergence of optimal distributed source controls]\label{thm:main-control-rate}
Assume data \eqref{eq:ocp-prepared}, a fixed target $y_d\in\YY$, a fixed $\gamma>0$, and a nonempty closed convex set $\UU_{\ad}\subset L^2H$. For each $(c_0,\delta)\in[0,\bar c_0]\times[0,\bar\delta]$ let $\bar m_{p,c_0,\delta}$ be the unique optimal distributed source control for problem \eqref{eq:cost} in the setting \eqref{eq:source-family}. Then
\begin{equation}
\norm{\bar m_{p,c_0,\delta}-\bar m_{p,0,0}}_{L^2H}
\le C_\gamma(c_0+\delta).
\label{eq:main-control-rate}
\end{equation}
The corresponding optimal states satisfy
\begin{equation}
\norm{(\bar u_{c_0,\delta},\bar p_{c_0,\delta})-(\bar u_{0,0},\bar p_{0,0})}_{\YY_V}
\le C_\gamma(c_0+\delta),
\label{eq:main-state-rate}
\end{equation}
and the optimal values satisfy
\begin{equation}
\abs{j_{c_0,\delta}(\bar m_{p,c_0,\delta})-j_{0,0}(\bar m_{p,0,0})}
\le C_\gamma(c_0+\delta).
\label{eq:value-rate}
\end{equation}
The constants may depend on the fixed data, the target, an arbitrary comparison element of $\UU_{\ad}$, $\gamma^{-1}$, and the parameter rectangle, but not on $c_0$ or $\delta$. The estimate is independent of active-set and complementarity assumptions.
\end{theorem}
\begin{proof}
The reduced objective is \eqref{eq:abstract-perturbed-cost} with $\theta=(c_0,\delta)$, $0=(0,0)$, and $y_\theta=y^b_{c_0,\delta}-y_d\in\YY$. Since $y_\theta-y_0=y^b_{c_0,\delta}-y^b_{0,0}$ and $\YY_V\hookrightarrow\YY$, Proposition~\ref{prop:operator-rate} gives
\begin{equation*}
\sup_{c_0,\delta}\norm{\mathcal S_{c_0,\delta}}_{\calL(\UU,\YY)}<\infty,
\quad
\norm{\mathcal S_{c_0,\delta}-\mathcal S_{0,0}}_{\calL(\UU,\YY)}+\norm{y_\theta-y_0}_\YY
\le C(c_0+\delta).
\end{equation*}
Apply Lemma~\ref{lem:perturbation}, with $m_0=\bar m_{p,0,0}$ fixed, to get \eqref{eq:main-control-rate}.

For the corresponding states, insert the intermediate state at control $\bar m_{p,0,0}$:
\begin{align*}
&\norm{\mathcal S_{c_0,\delta}\bar m_{p,c_0,\delta}+y^b_{c_0,\delta}
-\mathcal S_{0,0}\bar m_{p,0,0}-y^b_{0,0}}_{\YY_V}\\
&\ \le\norm{\mathcal S_{c_0,\delta}(\bar m_{p,c_0,\delta}-\bar m_{p,0,0})}_{\YY_V}
+\norm{(\mathcal S_{c_0,\delta}-\mathcal S_{0,0})\bar m_{p,0,0}}_{\YY_V}
+\norm{y^b_{c_0,\delta}-y^b_{0,0}}_{\YY_V}.
\end{align*}
The uniform operator bound, \eqref{eq:operator-rate}, and \eqref{eq:main-control-rate} give \eqref{eq:main-state-rate}.

To prove convergence of the optimal values, fix a comparison $m_{p,*}\in\UU_{\ad}$. Since $j_{c_0,\delta}(\bar m_{p,c_0,\delta})\le j_{c_0,\delta}(m_{p,*})$ and $j_{c_0,\delta}(m_{p,*})$ remains uniformly bounded by Proposition~\ref{prop:operator-rate}, the optimal controls are uniformly bounded in $\UU$. Hence the associated states and tracking residuals are uniformly bounded in $\YY$. The difference of squared norms satisfies
\begin{equation*}
\abs{\norm x_\YY^2-\norm y_\YY^2}
\le(\norm x_\YY+\norm y_\YY)\norm{x-y}_\YY.
\end{equation*}
Apply this inequality to the two optimal residuals and use \eqref{eq:main-state-rate}. For the regularization term use
\begin{equation*}
\abs{\norm{\bar m_{p,c_0,\delta}}_\UU^2-\norm{\bar m_{p,0,0}}_\UU^2}
\le(\norm{\bar m_{p,c_0,\delta}}_\UU+\norm{\bar m_{p,0,0}}_\UU)
\norm{\bar m_{p,c_0,\delta}-\bar m_{p,0,0}}_\UU.
\end{equation*}
This proves \eqref{eq:value-rate}.
\end{proof}

\begin{remark}[Energy and tracking norms]\label{rem:newrate}
Using only \eqref{eq:combined-energy} in Proposition~\ref{prop:operator-rate} would give the weaker estimate $O(c_0+\sqrt\delta)$ for the state maps and controls. The $O(c_0+\delta)$ rate in Theorem~\ref{thm:main-control-rate} depends crucially on the duality identity \eqref{eq:delta-duality-exact} and the $\delta$-uniform bound \eqref{eq:uniform-visc-velocity}.
\end{remark}

\begin{remark}[Observation operators]
The first-order estimate also holds for tracking through a bounded linear observation operator $\mathcal O:\YY_V\to Z$, with a target in $Z$. This includes localized displacement and pressure observations and the elastic stress contribution $\mathcal C_{\rm e}\eps(u)-\alpha pI$. For pressure-gradient observations, \eqref{eq:combined-energy} gives only $O(c_0+\sqrt\delta)$ for the observed state, and the exponent $1/2$ is sharp by Proposition~\ref{prop:source-sharp}.
\end{remark}

\subsection{Extension to boundary controls}\label{sec:boundary-controls}
Let $\mathrm{Tr}_u:V\to L^2(\Gamma_{u,N})^d$ denote the trace operator, and let $\mathrm{Tr}_p$ be as in Lemma~\ref{lem:mult-trace}. By \eqref{eq:loads}, a traction $t_N$ and an outward flux $\psi_N$ generate the loads $\mathrm{Tr}_u^*t_N\in V^*$ and $-\mathrm{Tr}_p^*\psi_N\in Q^*$, where $\dual{\mathrm{Tr}_u^*t}v=(t,\mathrm{Tr}_uv)_{\Gamma_{u,N}}$ for $v\in V$ and $\dual{\mathrm{Tr}_p^*\psi}q=(\psi,\mathrm{Tr}_pq)_{\Gamma_{p,N}}$ for $q\in Q$. Throughout this subsection the background data satisfy \eqref{eq:ocp-prepared}, and $y^b_{c_0,\delta}$ is as in Section~\ref{subsec:source}.

\paragraph{Traction controls}
To compare traction controls down to $\delta=0$, we choose a control space for which the induced mechanical load belongs to $H^1(0,T;V^*)$, as required by the state theory at $\delta=0$. We prescribe the initial traction $t_N(0)=0$, so that the prepared displacement is independent of the control and the endpoint term \eqref{eq:initial-correction} vanishes for every control variation. Mechanical controls that are only $L^2$ in time, and control-dependent initial equilibria, are not treated here. Let
\begin{equation}
\begin{aligned}
&\UU_{t_N}:=\bigl\{t_N\in H^1(0,T;L^2(\Gamma_{u,N})^d):\ t_N(0)=0\bigr\},\\
&\norm{t_N}_{\UU_{t_N}}:=\norm{\partial_tt_N}_{L^2(0,T;L^2(\Gamma_{u,N})^d)},\\
&F=F_b+\mathrm{Tr}_u^*t_N,\qquad G=G_b,
\end{aligned}
\label{eq:traction-family}
\end{equation}
and let $\mathcal S^{t_N}_{c_0,\delta}:\UU_{t_N}\to\YY_V$ be the state generated by $F=\mathrm{Tr}_u^*t_N$, with $G=0$ and homogeneous initial data. Since $t_N(0)=0$, the norm $\norm\cdot_{\UU_{t_N}}$ is equivalent to the $H^1(0,T;L^2(\Gamma_{u,N})^d)$ norm on $\UU_{t_N}$, and the full load in \eqref{eq:traction-family} satisfies $F(0)=F_b(0)$, which gives the independence of the prepared displacement in \eqref{eq:ocp-prepared}.

\begin{corollary}[Convergence of optimal traction controls]\label{cor:traction}
Assume data \eqref{eq:ocp-prepared}, a fixed target $y_d\in\YY$, a fixed $\gamma>0$, and a nonempty closed convex set $\UU_{\ad}\subset\UU_{t_N}$. Then, uniformly for $(c_0,\delta)\in[0,\bar c_0]\times[0,\bar\delta]$,
\begin{equation}
\norm{\mathcal S^{t_N}_{c_0,\delta}}_{\calL(\UU_{t_N},\YY_V)}\le C,\qquad
\norm{\mathcal S^{t_N}_{c_0,\delta}-\mathcal S^{t_N}_{0,0}}_{\calL(\UU_{t_N},\YY_V)}\le C(c_0+\delta),
\label{eq:traction-operator-rate}
\end{equation}
and the optimal traction controls $\bar t_{N,c_0,\delta}$ for \eqref{eq:cost} with $\UU=\UU_{t_N}$ and \eqref{eq:traction-family}, with optimal states $(\bar u_{c_0,\delta},\bar p_{c_0,\delta})$, satisfy
\begin{equation}
\begin{aligned}
&\norm{\bar t_{N,c_0,\delta}-\bar t_{N,0,0}}_{\UU_{t_N}}
+\norm{(\bar u_{c_0,\delta},\bar p_{c_0,\delta})-(\bar u_{0,0},\bar p_{0,0})}_{\YY_V}\\
&\qquad+\abs{j_{c_0,\delta}(\bar t_{N,c_0,\delta})-j_{0,0}(\bar t_{N,0,0})}
\le C_\gamma(c_0+\delta),
\end{aligned}
\label{eq:traction-control-rate}
\end{equation}
with constants as in Theorem~\ref{thm:main-control-rate}. In the corresponding optimality systems, the endpoint term \eqref{eq:initial-correction} vanishes.
\end{corollary}
\begin{proof}
For $t_N\in\UU_{t_N}$, the data $F=\mathrm{Tr}_u^*t_N$, $G=0$, $p_0=0$ satisfy Assumption~\ref{ass:regular-data}: $F\in H^1V^*$ and $F(0)=0$, so the prepared displacement is $u_0=0$, and $\mathcal D=\norm{F_t}_{L^2V^*}\le\norm{\mathrm{Tr}_u}\norm{t_N}_{\UU_{t_N}}$. Corollary~\ref{cor:combined-rate} gives the second bound in \eqref{eq:traction-operator-rate}. For the elastic state, \eqref{eq:pressure-state-est} with $p_0=0$ gives $\norm p_{L^2Q}\le C\norm{F_t}_{L^2V^*}$; since $F(0)=0$, $\norm F_{L^2V^*}\le T\norm{F_t}_{L^2V^*}$, and $u=A_{\ee}^{-1}(F+B^*p)$ yields $\norm{\mathcal S^{t_N}_{0,0}}\le C$. The first bound follows by the triangle inequality. The background term is the same as in Proposition~\ref{prop:operator-rate}, so the proof of Theorem~\ref{thm:main-control-rate} applies verbatim with $\UU=\UU_{t_N}$. Finally, $F_\eta(0)=\mathrm{Tr}_u^*\eta(0)=0$ for every variation $\eta\in\UU_{t_N}$, so \eqref{eq:initial-correction} vanishes.
\end{proof}

\paragraph{Flux controls}
Let
\begin{equation}
\UU_{\psi_N}:=L^2(0,T;L^2(\Gamma_{p,N})),\qquad F=F_b,\qquad G=G_b-\mathrm{Tr}_p^*\psi_N,
\label{eq:flux-family}
\end{equation}
and let $\mathcal S^{\psi_N}_{c_0,\delta}:\UU_{\psi_N}\to\YY_V$ be the state generated by $G_{\psi_N}:=-\mathrm{Tr}_p^*\psi_N$, with $F=0$ and homogeneous initial data. Since $G_{\psi_N}\in L^2Q^*$ but in general $G_{\psi_N}\notin L^2H$, Lemma~\ref{lem:uniform-reg} does not apply to $G_{\psi_N}$. Instead we use the pressure-adjoint estimates of Theorem~\ref{thm:adjoint-limit}.

\begin{theorem}[Convergence of the flux control-to-state operators]\label{thm:flux-rate}
Uniformly for $(c_0,\delta)\in[0,\bar c_0]\times[0,\bar\delta]$,
\begin{equation}
\norm{\mathcal S^{\psi_N}_{c_0,\delta}}_{\calL(\UU_{\psi_N},\YY_V)}\le C,\qquad
\norm{\mathcal S^{\psi_N}_{c_0,\delta}-\mathcal S^{\psi_N}_{0,0}}_{\calL(\UU_{\psi_N},\YY_V)}\le C\bigl(c_0+\delta^{3/4}\bigr).
\label{eq:flux-operator-rate}
\end{equation}
\end{theorem}
\begin{proof}
Since $\mathcal S^{\psi_N}_{c_0,\delta}\psi_N=\mathcal T_{c_0,\delta}G_{\psi_N}$ and $\norm{G_{\psi_N}}_{L^2Q^*}\le\norm{\mathrm{Tr}_p}\norm{\psi_N}_{\UU_{\psi_N}}$, Lemma~\ref{lem:Qstar-stability} gives the first bound in \eqref{eq:flux-operator-rate}. For the second, let $R_u\in L^2V^*$ and $R_p\in L^2H$, and let $\tilde p_{c_0,\delta}$ be the corresponding pressure adjoints. By \eqref{eq:adj-duality-G} with $G=G_{\psi_N}$,
\begin{equation*}
\int_0^T\!\left\{\dual{R_u}w+(R_p,r)_H\right\}\dd t
=-\int_0^T\bigl(\psi_N,\mathrm{Tr}_p(\tilde p_{c_0,\delta}-\tilde p_{0,0})\bigr)_{\Gamma_{p,N}}\dd t
\end{equation*}
for $(w,r)=\bigl(\mathcal S^{\psi_N}_{c_0,\delta}-\mathcal S^{\psi_N}_{0,0}\bigr)\psi_N$.
By \eqref{eq:adj-joint}, the right-hand side is bounded by $C(c_0+\delta^{3/4})\bigl(\norm{R_u}_{L^2V^*}+\norm{R_p}_{L^2H}\bigr)\norm{\psi_N}_{\UU_{\psi_N}}$. Since $L^2V^*\times L^2H$ is the dual of $\YY_V$, taking the supremum over $(R_u,R_p)$ in its unit ball proves the second bound in \eqref{eq:flux-operator-rate}.
\end{proof}

\begin{corollary}[Convergence of optimal flux controls]\label{cor:flux}
Assume data \eqref{eq:ocp-prepared}, a fixed target $y_d\in\YY$, a fixed $\gamma>0$, and a nonempty closed convex set $\UU_{\ad}\subset\UU_{\psi_N}$. Then the optimal flux controls $\bar\psi_{N,c_0,\delta}$ for \eqref{eq:cost} with $\UU=\UU_{\psi_N}$ and \eqref{eq:flux-family}, with optimal states $(\bar u_{c_0,\delta},\bar p_{c_0,\delta})$, satisfy
\begin{equation}
\begin{aligned}
&\norm{\bar\psi_{N,c_0,\delta}-\bar\psi_{N,0,0}}_{\UU_{\psi_N}}
+\norm{(\bar u_{c_0,\delta},\bar p_{c_0,\delta})-(\bar u_{0,0},\bar p_{0,0})}_{\YY_V}\\
&\qquad+\abs{j_{c_0,\delta}(\bar\psi_{N,c_0,\delta})-j_{0,0}(\bar\psi_{N,0,0})}
\le C_\gamma\bigl(c_0+\delta^{3/4}\bigr),
\end{aligned}
\label{eq:flux-control-rate}
\end{equation}
with constants as in Theorem~\ref{thm:main-control-rate}.
\end{corollary}
\begin{proof}
Theorem~\ref{thm:flux-rate} and the background estimate of Proposition~\ref{prop:operator-rate}, which is $O(c_0+\delta)$ and hence $O(c_0+\delta^{3/4})$ on the parameter rectangle, give \eqref{eq:convex-perturb-assumption} with $\omega_{c_0,\delta}=C(c_0+\delta^{3/4})$. The proof of Theorem~\ref{thm:main-control-rate} then applies verbatim.
\end{proof}

\begin{proposition}[Sharpness of the flux exponent]\label{prop:flux-sharp}
For the configuration constructed in Appendix~\ref{subsec:example} for Proposition~\ref{prop:source-sharp}, there is a constant $\kappa>0$ such that for every sufficiently small $\delta>0$ there is a flux $\psi_N$ of unit norm in $\UU_{\psi_N}$ for which
\begin{equation*}
\norm{\bigl(\mathcal S^{\psi_N}_{0,\delta}-\mathcal S^{\psi_N}_{0,0}\bigr)\psi_N}_{\YY_V}
\ge\norm{\bigl(\mathcal S^{\psi_N}_{0,\delta}-\mathcal S^{\psi_N}_{0,0}\bigr)\psi_N}_{\YY}
\ge \kappa\delta^{3/4}.
\end{equation*}
Hence the exponent $3/4$ in \eqref{eq:flux-operator-rate} cannot be improved in general.
\end{proposition}

For the unconstrained problem, the same exponent is sharp in the operator norm of the map from the target to the optimal flux.

\begin{corollary}[Sharpness for the target-to-optimal-flux map]\label{cor:flux-sharp-control}
In the configuration of Proposition~\ref{prop:flux-sharp}, let $F_b=0$, $G_b=0$, $p_0=0$, and $\UU_{\ad}=\UU_{\psi_N}$, and fix $\gamma>0$. For $0\le\delta\le\bar\delta$ let $\Psi_\delta:\YY\to\UU_{\psi_N}$ assign to a target $y_d\in\YY$ the optimal flux control $\bar\psi_{N,0,\delta}$ for \eqref{eq:cost} with $\UU=\UU_{\psi_N}$ and \eqref{eq:flux-family}. Then $\Psi_\delta$ is linear and bounded, and there are constants $C_\gamma,c_\gamma,\delta_\gamma>0$, depending on $\gamma$, such that
\begin{equation}
\begin{aligned}
\norm{\Psi_\delta-\Psi_0}_{\calL(\YY,\UU_{\psi_N})}&\le C_\gamma\delta^{3/4}, &&0<\delta\le\bar\delta,\\
\norm{\Psi_\delta-\Psi_0}_{\calL(\YY,\UU_{\psi_N})}&\ge c_\gamma\delta^{3/4}, &&0<\delta\le\min\{\delta_\gamma,T,\bar\delta\}.
\end{aligned}
\label{eq:flux-control-sharp}
\end{equation}
\end{corollary}
The proofs of Proposition~\ref{prop:flux-sharp} and Corollary~\ref{cor:flux-sharp-control} are given in Appendix~\ref{app:flux-sharpness}.

\begin{remark}[Optimal adjoints]
Theorem~\ref{thm:adjoint-limit} compares adjoints with the same residuals. For the optimal adjoints the residuals $\bar R_{c_0,\delta}:=(\bar u_{c_0,\delta}-u_d,\bar p_{c_0,\delta}-p_d)$ depend on the parameters. Write $\tilde p_{c_0,\delta}[R]$ for the pressure adjoint with residual $R=(R_u,R_p)$. By linearity,
\begin{equation*}
\tilde p_{c_0,\delta}[\bar R_{c_0,\delta}]-\tilde p_{0,0}[\bar R_{0,0}]
=\tilde p_{c_0,\delta}[\bar R_{c_0,\delta}-\bar R_{0,0}]
+\bigl(\tilde p_{c_0,\delta}-\tilde p_{0,0}\bigr)[\bar R_{0,0}].
\end{equation*}
By \eqref{eq:adj-uniform-Q}, the first term is bounded in $L^2Q$, hence in $L^2H$, by $C$ times the norm of $\bar R_{c_0,\delta}-\bar R_{0,0}$ in $L^2V^*\times L^2H$, which is at most $C$ times the optimal-state error in $\YY_V$. By \eqref{eq:adj-joint}, the second term is bounded in $L^2H$ by $C(c_0+\delta)$ times the norm of $\bar R_{0,0}$. Hence the optimal pressure adjoints converge in $L^2H$ at the rate $O(c_0+\delta)$ for distributed source and traction controls, by \eqref{eq:main-state-rate} and \eqref{eq:traction-control-rate}, and at the rate $O(c_0+\delta^{3/4})$ for flux controls, by \eqref{eq:flux-control-rate}.
\end{remark}

\section{Conclusion}
We have analyzed a linear quasi-static poroviscoelastic optimal control problem for all $c_0,\delta\geq0$. The choice of evolution variable for the state and adjoint equations depends on the parameter regime. In particular, for $c_0=0$ and $\delta>0$ the adjoint terminal condition is imposed on $\delta A_{\vv}\tilde u+B^*\tilde p$, whereas for $\delta=0$ the pressure evolves with the mass operator $M_{c_0}=c_0I+BA_{\ee}^{-1}B^*$.

Under the standing assumption on $\Gamma_{u,N}$, which yields the pressure-to-dilation inf--sup condition, and for data satisfying Assumption~\ref{ass:regular-data}, the states converge at the rate $O(c_0+\sqrt\delta)$ in the energy norm and $O(c_0+\delta)$ in $L^2V\times L^2H$, and both exponents of $\delta$ are sharp. The second estimate yields operator-norm convergence of the state maps and the same first-order rate for constrained distributed source controls with positive quadratic regularization. The same rate holds for traction controls in an $H^1$-in-time space with zero initial trace. For fixed observation residuals, the adjoint pair is bounded uniformly in both parameters, the mechanical adjoint converges in $L^2(0,T;V)$, and the pressure adjoint converges at the rate $O(c_0+\delta)$ in $L^2(0,T;L^2(\Omega))$ and $O(c_0+\sqrt\delta)$ in $L^2(0,T;Q)$. A multiplicative trace inequality turns these estimates into the rate $O(c_0+\delta^{3/4})$ for normal Darcy flux controls. An explicit example shows that this exponent cannot be improved for the flux-to-state operators, nor, without control constraints, for the map from the target to the optimal flux control. We do not establish a matching lower bound for optimal flux controls under active pointwise constraints, or for isotropic elasticity tensors.

\appendix
\section{Proofs of the sharpness results}\label{app:sharpness}
\subsection{The displacement and the mechanical adjoint}\label{app:mechanical-sharpness}

\begin{proof}[Proof of Theorem~\ref{thm:visc-sharp}]
Choose $w\in V\setminus\{0\}$ with $Bw=0$, for instance a nonzero divergence-free field in $H^1_0(\Omega)^d\subset V$. Set
\begin{equation*}
f_\delta(t):=\delta^{-1/2}(t-T+\delta)_+,\qquad
F:=f_\delta\,A_{\ee}w,\qquad G:=0,\qquad p_0:=0.
\end{equation*}
Since $\delta\le T$, $f_\delta(0)=0$, $f_\delta\in H^1(0,T)$, $\partial_tf_\delta=\delta^{-1/2}\chi_{(T-\delta,T)}$, and $\norm{\partial_tf_\delta}_{L^2(0,T)}=1$. Thus $F\in H^1V^*$, the prepared displacement \eqref{eq:prepared} is $u_0=0$, and
\begin{equation*}
\mathcal D=\norm{F_t}_{L^2V^*}=\norm{A_{\ee}w}_{V^*},
\qquad
\underline C_{\ee}\norm w_V\le\mathcal D\le\overline C_{\ee}\norm w_V.
\end{equation*}
The pair $(f_\delta w,0)$ satisfies $A_{\ee}(f_\delta w)=F$, $(c_0\cdot0+B(f_\delta w))_t=\partial_tf_\delta\,Bw=0$, and $p(0)=0$; by uniqueness in Theorem~\ref{thm:state-storage}, $(u_{c_0,0},p_{c_0,0})=(f_\delta w,0)$. Let $\nu$ solve $\delta\nu_t+\nu=f_\delta$, $\nu(0)=0$. Since $\mathcal C_{\rm v}=\mathcal C_{\rm e}$, we have $A_{\vv}=A_{\ee}$, and the pair $(\nu w,0)$ satisfies $\delta A_{\vv}(\nu w)_t+A_{\ee}(\nu w)=(\delta\nu_t+\nu)A_{\ee}w=F$, $(c_0\cdot0+B(\nu w))_t=0$, $\nu(0)w=0=u_0$, and $p(0)=0=p_0$ when $c_0>0$. By uniqueness in Theorems~\ref{thm:state-full} and \ref{thm:state-visc}, $(u_{c_0,\delta},p_{c_0,\delta})=(\nu w,0)$.

Hence $e_p=0$ and $e_u=b_\delta w$, where $b_\delta:=\nu-f_\delta$ solves $\delta b_\delta'+b_\delta=-\delta\,\partial_tf_\delta$, $b_\delta(0)=0$. Therefore $b_\delta=0$ on $[0,T-\delta]$, and, with $s=t-T+\delta\in[0,\delta]$,
\begin{equation*}
b_\delta(t)=-\sqrt\delta\bigl(1-\exp(-s/\delta)\bigr),\qquad
b_\delta'(t)=-\delta^{-1/2}\exp(-s/\delta).
\end{equation*}
It follows that $\norm{e_u}_{L^\infty V}\ge\abs{b_\delta(T)}\norm w_V=(1-\exp(-1))\sqrt\delta\norm w_V$ and
\begin{equation*}
\norm{\partial_t e_{u}}_{L^2V}^2
=\norm w_V^2\,\delta^{-1}\int_0^\delta\exp(-2s/\delta)\dd s
=\tfrac12\bigl(1-\exp(-2)\bigr)\norm w_V^2.
\end{equation*}
Together with $\mathcal D\le\overline C_{\ee}\norm w_V$, this proves \eqref{eq:visc-sharp}. If $\norm{e_u}_{L^\infty V}\le C\delta^\sigma\mathcal D$ held for all such data, then $(1-\exp(-1))\delta^{1/2}\le C\overline C_{\ee}\delta^\sigma$ for all small $\delta>0$, which is impossible for $\sigma>1/2$; the same argument applies to $\sqrt\delta\norm{\partial_t e_{u}}_{L^2V}$.

The same example also illustrates the dependence on the time norm. Since $b_\delta$ vanishes on $[0,T-\delta]$ and $\abs{b_\delta}\le\sqrt\delta$, we have $\norm{e_u}_{L^2V}\le\delta\norm w_V\le\delta\mathcal D/\underline C_{\ee}$, in agreement with Theorem~\ref{thm:visc-dual-rate}.
\end{proof}

\begin{proof}[Sharpness in Corollary~\ref{cor:mech-adjoint-rate}]
Let $c_0=0$ and $\mathcal C_{\rm v}=\mathcal C_{\rm e}$, choose $v\in V\setminus\{0\}$ with $Bv=0$ as in the proof of Theorem~\ref{thm:visc-sharp}, and take $R_u=A_{\ee}v$ and $R_p=0$. Then the right-hand side of \eqref{eq:adj-reduced-pressure} is $Bv=0$, so $\tilde p_{0,0}=0$ by Theorem~\ref{thm:adj-elastic} and $\tilde u_{0,0}=v$ by \eqref{eq:adj-elim-m}. For $\delta>0$, the pair $\tilde u=a_\delta(t)v$, $\tilde p=0$ with $a_\delta(t):=1-\exp(-(T-t)/\delta)$ satisfies \eqref{eq:visc-adjoint-form}, because $a_\delta-\delta a_\delta'=1$ and $Bv=0$, and $z=\delta a_\delta A_{\ee}v\in H^1V^*$ with $z(T)=0$; by Theorem~\ref{thm:adj-visc} it is the adjoint. Hence
\begin{equation*}
\norm{\tilde u_{0,\delta}-\tilde u_{0,0}}_{L^2V}^2=\norm v_V^2\int_0^T\exp\bigl(-2(T-t)/\delta\bigr)\dd t
=\frac\delta2\bigl(1-\exp(-2T/\delta)\bigr)\norm v_V^2,
\end{equation*}
which is at least $\frac{1-e^{-2}}{2T}\,\delta\,\norm{\tilde u_{0,0}}_{H^1V}^2$ for $0<\delta\le T$. Thus $\tilde u_{0,\delta}(T)=a_\delta(T)v=0$, whereas $\tilde u_{0,0}(T)=v$.
\end{proof}

\subsection{A two-dimensional example for pressure and flux}\label{subsec:example}

We construct exact solutions in a two-dimensional configuration that satisfies all assumptions of Section~\ref{sec:setting}. Let
\begin{equation*}
\begin{gathered}
\Omega=(0,1)^2,\qquad \Gamma_{u,D}=\{0\}\times(0,1),\qquad \Gamma_{p,D}=\{1\}\times(0,1),\\
\alpha=1,\qquad \mathbf K=I,\qquad c_0=0,
\end{gathered}
\end{equation*}
so that $\Gamma_{u,N}$ contains the relatively open set $\{1\}\times(0,1)$. For symmetric $E\in\mathbb R^{2\times2}$ let
\begin{equation*}
\mathcal C_{\rm e}E=\mathcal C_{\rm v}E:=\begin{pmatrix}E_{11}+E_{22}&2E_{12}\\ 2E_{12}&E_{11}+2E_{22}\end{pmatrix},
\qquad
\mathcal C_{\rm e}E:E=(E_{11}+E_{22})^2+E_{22}^2+4E_{12}^2.
\end{equation*}
This tensor is symmetric and positive definite on symmetric matrices, and $\mathcal C_{\rm e}E=sI$ for $E=\mathrm{diag}(s,0)$. Let $\phi_n(x)=\sqrt2\cos((n+\frac12)\pi x)$ and $\rho_n=(n+\frac12)^2\pi^2$, $n\ge0$; these are the orthonormal eigenfunctions and eigenvalues of $-\partial_x^2$ on $(0,1)$ with $\phi_n'(0)=0$ and $\phi_n(1)=0$, and $\phi_n(0)=\sqrt2$.

\begin{lemma}[Reduction to one dimension]\label{lem:example-reduction}
Let $\delta\ge0$, $g\in L^2(0,T;L^2(0,1))$, and $\psi\in L^2(0,T)$, and consider $F=0$, zero initial data, and the fluid load
\begin{equation*}
\dual Gq=\int_\Omega g(x,t)\,q(x,y)\dd x\dd y-\psi(t)\int_0^1q(0,y)\dd y,\qquad q\in Q,
\end{equation*}
that is, the source $g$ and the outward flux $\psi_N=\psi$ on $\{0\}\times(0,1)$ and $\psi_N=0$ on the rest of $\Gamma_{p,N}$. Put $G_n:=(g,\phi_n)_{L^2(0,1)}-\sqrt2\,\psi$ and let $U_n$ solve
\begin{equation*}
(1+\delta\rho_n)U_n'+\rho_nU_n=G_n,\qquad U_n(0)=0,
\qquad\text{and}\qquad
P_n:=\frac{U_n+\delta G_n}{1+\delta\rho_n}.
\end{equation*}
Then the state $(u_{0,\delta},p_{0,\delta})$ of Theorem~\ref{thm:state-visc} (for $\delta>0$) or Theorem~\ref{thm:state-storage} (for $\delta=0$) is
\begin{equation*}
\begin{gathered}
u_{0,\delta}(x,y,t)=\Bigl(\int_0^xU(s,t)\dd s,\ 0\Bigr),\qquad p_{0,\delta}(x,y,t)=P(x,t),\\
U=\sum_nU_n\phi_n,\qquad P=\sum_nP_n\phi_n.
\end{gathered}
\end{equation*}
\end{lemma}
\begin{proof}
Let $Q_{\rm 1D}:=\{\vartheta\in H^1(0,1):\vartheta(1)=0\}$. Multiplying the modal equation by $U_n$ and integrating over $(0,T)$ with $U_n(0)=0$ gives
\begin{equation*}
\frac{1+\delta\rho_n}2\abs{U_n(T)}^2+\rho_n\norm{U_n}_{L^2(0,T)}^2=(G_n,U_n)_{L^2(0,T)},
\end{equation*}
so $\rho_n\norm{U_n}_{L^2(0,T)}\le\norm{G_n}_{L^2(0,T)}$; since $U_n'=(G_n-\rho_nU_n)/(1+\delta\rho_n)$ and $P_n=(U_n+\delta G_n)/(1+\delta\rho_n)$, this gives
\begin{equation*}
\norm{U_n'}_{L^2(0,T)}\le\frac{2\norm{G_n}_{L^2(0,T)}}{1+\delta\rho_n},\qquad
\norm{P_n}_{L^2(0,T)}\le\frac{\norm{G_n}_{L^2(0,T)}}{\rho_n}.
\end{equation*}
Since $\sum_n\rho_n^{-1}\norm{G_n}_{L^2(0,T)}^2\le C(\norm g_{L^2L^2}^2+\norm\psi_{L^2(0,T)}^2)$, we obtain $P\in L^2(0,T;Q_{\rm 1D})$ and $U_t=\sum_nU_n'\phi_n\in L^2(0,T;Q_{\rm 1D}^*)$ for every $\delta\ge0$; for $\delta>0$, moreover, $\sum_n\norm{U_n'}_{L^2(0,T)}^2<\infty$ because $\rho_n\sim n^2$, so $U\in H^1(0,T;L^2(0,1))$. For $\delta=0$ one has $U=P\in L^2(0,T;Q_{\rm 1D})\cap H^1(0,T;Q_{\rm 1D}^*)\subset C([0,T];L^2(0,1))$. In all cases $U(0)=0$, and the modal equations are equivalent to
\begin{subequations}\label{eq:example-1d}
\begin{align}
\delta U_t+U&=P,
\label{eq:example-1d-mech}\\
\dual{U_t}{\vartheta}_{Q_{\rm 1D}^*,Q_{\rm 1D}}+(P_x,\vartheta')_{L^2(0,1)}&=(g,\vartheta)_{L^2(0,1)}-\psi\,\vartheta(0)
\quad\text{for all }\vartheta\in Q_{\rm 1D}.
\label{eq:example-1d-fluid}
\end{align}
\end{subequations}
With $u$ as above, $\eps(u)=\mathrm{diag}(U,0)$, $\divergence u=U$, and hence, by \eqref{eq:example-1d-mech}, $\mathcal C_{\rm e}\eps(u)+\delta\mathcal C_{\rm v}\eps(u_t)-pI=(U+\delta U_t-P)I=0$. Thus $\delta a_{\vv}(u_t,v)+a_{\ee}(u,v)-(p,Bv)_H=\int_\Omega(U+\delta U_t-P)\divergence v\dd x\dd y=0$ for every $v\in V$, which is \eqref{eq:state-mech} with $F=0$; moreover $u=0$ on $\Gamma_{u,D}$. For $q\in Q$ put $\vartheta(x):=\int_0^1q(x,y)\dd y$. Then $\vartheta\in Q_{\rm 1D}$, $\vartheta(0)=\int_0^1q(0,y)\dd y$, $(P_x,\vartheta')_{L^2(0,1)}=a_p(p,q)$, $(g,\vartheta)_{L^2(0,1)}=\int_\Omega gq\dd x\dd y$, and $\dual{(Bu)_t}q=\dual{U_t}\vartheta_{Q_{\rm 1D}^*,Q_{\rm 1D}}$ because $Bu=U$ does not depend on $y$. Hence \eqref{eq:example-1d-fluid} gives \eqref{eq:state-fluid} with $c_0=0$. The pair $(u,p)$ lies in the solution class of Theorem~\ref{thm:state-visc} or \ref{thm:state-storage}, and uniqueness in these theorems identifies it with the state.
\end{proof}

We write $P_n^{(\delta)}$ for the coefficient $P_n$ at viscosity $\delta\ge0$. For constant $G_n$ on $(0,\delta)$ the modal solutions are explicit: with $\lambda_n=\rho_n/(1+\delta\rho_n)$,
\begin{equation}
P_n^{(\delta)}(t)=\frac{G_n}{\rho_n}\Bigl(1-\frac{\exp(-\lambda_nt)}{1+\delta\rho_n}\Bigr),\qquad
P_n^{(0)}(t)=\frac{G_n}{\rho_n}\bigl(1-\exp(-\rho_nt)\bigr),\qquad 0<t<\delta.
\label{eq:example-modes}
\end{equation}
If $1\le\delta\rho_n\le2$ and $0<t<\delta/8$, then $\rho_nt\le\frac14$ and $(1+\delta\rho_n)^{-1}\le\frac12$, so subtracting the two formulas in \eqref{eq:example-modes} gives
\begin{equation}
\abs{P_n^{(\delta)}(t)-P_n^{(0)}(t)}\ge c_\sharp\frac{\abs{G_n}}{\rho_n},\qquad c_\sharp:=\exp(-1/4)-\tfrac12>0.
\label{eq:example-mode-gap}
\end{equation}
Since $\rho_{n+1}\le2\rho_n$ for $n\ge2$, for every $0<\delta\le\rho_2^{-1}$ there is $n$ with $1\le\delta\rho_n\le2$.

\begin{proof}[Proof of Proposition~\ref{prop:source-sharp}]
Let $0<\delta\le\min\{\rho_2^{-1},T,\bar\delta\}$ and choose $n$ so that $1\le\delta\rho_n\le2$. Set $F=0$, $p_0=0$, $\psi=0$, and $G_\delta(x,y,t)=\delta^{-1/2}\chi_{(0,\delta)}(t)\phi_n(x)$. Since $\phi_n$ has unit $L^2(0,1)$ norm, Assumption~\ref{ass:regular-data} holds with $\mathcal D=\norm{G_\delta}_{L^2H}=1$.

Only the mode $n$ is excited, with $G_n=\delta^{-1/2}$ on $(0,\delta)$. By \eqref{eq:example-mode-gap} and $\rho_n\le2/\delta$, $\abs{P_n^{(\delta)}-P_n^{(0)}}\ge c_\sharp\delta^{1/2}/2$ on $(0,\delta/8)$. For a pressure $r(x,y)=\sum_m r_m\phi_m(x)$ independent of $y$, Parseval's identity gives $\norm r_H^2=\sum_m\abs{r_m}^2$ and $\norm r_Q^2=\sum_m\rho_m\abs{r_m}^2$. Hence $\norm{p_{0,\delta}-p_{0,0}}_{L^2H}^2\ge\frac\delta8\cdot\frac{c_\sharp^2\delta}4$, and, since only the mode $n$ is excited, $\norm{p_{0,\delta}-p_{0,0}}_{L^2Q}^2=\rho_n\norm{p_{0,\delta}-p_{0,0}}_{L^2H}^2$ with $\rho_n\ge1/\delta$. Both bounds therefore hold with $c=c_\sharp/(4\sqrt2)$.
\end{proof}

\subsection{Boundary flux controls}\label{app:flux-sharpness}

\begin{proof}[Proof of Proposition~\ref{prop:flux-sharp}]
Take $0<\delta\le\min\{\delta_*,T,\bar\delta\}$, where $\delta_*:=\bigl((\sqrt2-1)/(2\pi)\bigr)^2$. In the configuration of Appendix~\ref{subsec:example}, take $g=0$ and choose the outward flux $\psi_N=-\delta^{-1/2}\chi_{(0,\delta)}(t)$ on $\{0\}\times(0,1)$ and zero on the rest of $\Gamma_{p,N}$. Then $\norm{\psi_N}_{\UU_{\psi_N}}=1$, and Lemma~\ref{lem:example-reduction} gives $G_n=\sqrt2\,\delta^{-1/2}$ on $(0,\delta)$ for every $n$. For each $n$ with $1\le\delta\rho_n\le2$, \eqref{eq:example-mode-gap} gives $\abs{P_n^{(\delta)}-P_n^{(0)}}\ge c_\sharp\delta^{1/2}/\sqrt2$ on $(0,\delta/8)$. The condition $1\le\delta\rho_n\le2$ means that $n+\frac12$ lies in an interval of length $(\sqrt2-1)/(\pi\sqrt\delta)$, which contains at least $(\sqrt2-1)/(2\pi\sqrt\delta)$ such half-integers when $\delta\le\delta_*$. By Parseval's identity,
\begin{equation*}
\norm{p_{0,\delta}-p_{0,0}}_{L^2H}^2
\ge\frac\delta8\cdot\frac{\sqrt2-1}{2\pi\sqrt\delta}\cdot\frac{c_\sharp^2\delta}2
=\frac{(\sqrt2-1)c_\sharp^2}{32\pi}\,\delta^{3/2}.
\end{equation*}
The pressure component of the error is bounded by its $\YY$ norm, and $\norm\cdot_\YY\le\norm\cdot_{\YY_V}$. This proves the assertion with $\kappa=c_\sharp((\sqrt2-1)/(32\pi))^{1/2}$.
\end{proof}

\begin{proof}[Proof of Corollary~\ref{cor:flux-sharp-control}]
Write $S_\delta:=\mathcal S^{\psi_N}_{0,\delta}$, regarded as an element of $\calL(\UU_{\psi_N},\YY)$, let $S_\delta^*$ be its Hilbert adjoint, and put $D:=S_\delta-S_0$ and $H_\delta:=\gamma I+S_\delta^*S_\delta$. By Theorem~\ref{thm:flux-rate} and $\norm\cdot_\YY\le\norm\cdot_{\YY_V}$, there are $C_{\mathcal S}$ and $C_D$ with $\norm{S_\delta}\le C_{\mathcal S}$ and $\norm D\le C_D\delta^{3/4}$; moreover $H_\delta\ge\gamma I$, so $\norm{H_\delta^{-1}}\le\gamma^{-1}$. Since $y^b_{0,\delta}=0$ and $\UU_{\ad}=\UU_{\psi_N}$, the variational inequality \eqref{eq:abstract-VI} reduces to $S_\delta^*(S_\delta\bar\psi-y_d)+\gamma\bar\psi=0$, so $\Psi_\delta=H_\delta^{-1}S_\delta^*$.

\emph{Upper bound.} Since $S_\delta^*S_\delta-S_0^*S_0=S_\delta^*D+D^*S_0$,
\begin{equation*}
\Psi_\delta-\Psi_0=H_\delta^{-1}D^*-H_\delta^{-1}\bigl(S_\delta^*D+D^*S_0\bigr)H_0^{-1}S_0^*,
\end{equation*}
and therefore $\norm{\Psi_\delta-\Psi_0}\le(\gamma^{-1}+2C_{\mathcal S}^2\gamma^{-2})C_D\delta^{3/4}$.

\emph{Lower bound.} We use $\norm{\Psi_\delta-\Psi_0}=\norm{\Psi_\delta^*-\Psi_0^*}$ and $\Psi_\delta^*=S_\delta H_\delta^{-1}$. For $\xi\in\UU_{\psi_N}$, the relation $H_0=H_\delta-S_\delta^*D-D^*S_0$ gives
\begin{equation*}
\begin{aligned}
(\Psi_\delta^*-\Psi_0^*)H_0\xi
&=S_\delta\xi-S_\delta H_\delta^{-1}\bigl(S_\delta^*D+D^*S_0\bigr)\xi-S_0\xi\\
&=\bigl(I-S_\delta H_\delta^{-1}S_\delta^*\bigr)D\xi-S_\delta H_\delta^{-1}D^*S_0\xi.
\end{aligned}
\end{equation*}
Since $(\gamma I+S_\delta S_\delta^*)S_\delta=S_\delta H_\delta$, we have $S_\delta H_\delta^{-1}=(\gamma I+S_\delta S_\delta^*)^{-1}S_\delta$, and hence $I-S_\delta H_\delta^{-1}S_\delta^*=\gamma(\gamma I+S_\delta S_\delta^*)^{-1}$. This operator is self-adjoint on $\YY$ with spectrum in $[\gamma/(\gamma+C_{\mathcal S}^2),1]$. Consequently
\begin{equation*}
\norm{(\Psi_\delta^*-\Psi_0^*)H_0\xi}_\YY
\ge\frac{\gamma}{\gamma+C_{\mathcal S}^2}\norm{D\xi}_\YY-\frac{C_{\mathcal S}C_D}{\gamma}\,\delta^{3/4}\norm{S_0\xi}_\YY.
\end{equation*}
Let $\xi_\delta$ be the flux $\psi_N$ constructed in the proof of Proposition~\ref{prop:flux-sharp}. Then $\norm{\xi_\delta}_{\UU_{\psi_N}}=1$ and $\norm{D\xi_\delta}_\YY\ge\kappa\delta^{3/4}$ for $0<\delta\le\min\{\delta_*,T,\bar\delta\}$. We claim that $\norm{S_0\xi_\delta}_\YY\le\sqrt\delta$. By Lemma~\ref{lem:example-reduction} at viscosity zero, with $\phi_n$ and $\rho_n=((n+\frac12)\pi)^2$ as in Appendix~\ref{subsec:example}, $g=0$, and $\psi=-\delta^{-1/2}\chi_{(0,\delta)}$, the state $S_0\xi_\delta=(u_{0,0},p_{0,0})$ has pressure $P=U=\sum_nU_n\phi_n$, where $U_n'+\rho_nU_n=G_n:=\sqrt2\,\delta^{-1/2}\chi_{(0,\delta)}$ and $U_n(0)=0$. Thus $U_n$ is the convolution of $\exp(-\rho_n\,\cdot)$ with $G_n$, and Young's inequality, with the kernel in $L^2(0,\infty)$ and $G_n\in L^1(0,T)$, gives
\begin{equation*}
\norm{U_n}_{L^2(0,T)}\le\norm{\exp(-\rho_n\,\cdot)}_{L^2(0,\infty)}\norm{G_n}_{L^1(0,T)}
=(2\rho_n)^{-1/2}\sqrt2\,\delta^{1/2}=(\delta/\rho_n)^{1/2}.
\end{equation*}
By Parseval's identity and $\sum_n\rho_n^{-1}=\frac12$, $\norm{p_{0,0}}_{L^2H}^2=\sum_n\norm{U_n}_{L^2(0,T)}^2\le\delta/2$. Since $\abs{\int_0^xU(s,t)\dd s}\le\norm{U(\cdot,t)}_{L^2(0,1)}$, also $\norm{u_{0,0}}_{L^2(0,T;L^2(\Omega)^2)}^2\le\delta/2$, which proves the claim. Finally, $\norm{H_0\xi_\delta}_{\UU_{\psi_N}}\le\gamma+C_{\mathcal S}^2$, and hence
\begin{equation*}
\norm{\Psi_\delta-\Psi_0}
\ge\frac{\norm{(\Psi_\delta^*-\Psi_0^*)H_0\xi_\delta}_\YY}{\norm{H_0\xi_\delta}_{\UU_{\psi_N}}}
\ge\frac1{\gamma+C_{\mathcal S}^2}\Bigl(\frac{\gamma\kappa}{\gamma+C_{\mathcal S}^2}\,\delta^{3/4}-\frac{C_{\mathcal S}C_D}{\gamma}\,\delta^{5/4}\Bigr).
\end{equation*}
For $\delta\le\delta_\gamma:=\min\{\delta_*,(\gamma^2\kappa/(2C_{\mathcal S}C_D(\gamma+C_{\mathcal S}^2)))^2\}$ the bracket is at least $\gamma\kappa\delta^{3/4}/(2(\gamma+C_{\mathcal S}^2))$, which gives \eqref{eq:flux-control-sharp} with $c_\gamma:=\gamma\kappa/(2(\gamma+C_{\mathcal S}^2)^2)$ and $C_\gamma:=(\gamma^{-1}+2C_{\mathcal S}^2\gamma^{-2})C_D$.
\end{proof}

\section*{Acknowledgments}
The author used Claude Opus 5.5 (Anthropic) and ChatGPT 6 Astra (OpenAI) to draft and edit text, to develop and check parts of the proofs, and to check references; no AI-generated code, data, or figures appear. The author verified all mathematical content and references and assumes responsibility for all content.


\begin{thebibliography}{99}\small
\bibitem{BLMW2017}
\textsc{T.~B{\ae}rland, J.~J. Lee, K.-A. Mardal, and R.~Winther}, \emph{Weakly imposed symmetry and robust preconditioners for {B}iot's consolidation model}, Comput. Methods Appl. Math., 17 (2017), pp.~377--396, \url{https://doi.org/10.1515/cmam-2017-0016}.

\bibitem{Biot1941}
\textsc{M.~A. Biot}, \emph{General theory of three-dimensional consolidation}, J. Appl. Phys., 12 (1941), pp.~155--164, \url{https://doi.org/10.1063/1.1712886}.

\bibitem{BGSV2019}
\textsc{L.~Bociu, G.~Guidoboni, R.~Sacco, and M.~Verri}, \emph{On the role of compressibility in poroviscoelastic models}, Math. Biosci. Eng., 16 (2019), pp.~6167--6208, \url{https://doi.org/10.3934/mbe.2019308}.

\bibitem{BGSSW2016}
\textsc{L.~Bociu, G.~Guidoboni, R.~Sacco, and J.~T. Webster}, \emph{Analysis of nonlinear poro-elastic and poro-visco-elastic models}, Arch. Ration. Mech. Anal., 222 (2016), pp.~1445--1519, \url{https://doi.org/10.1007/s00205-016-1024-9}.

\bibitem{BMSS2024}
\textsc{L.~Bociu, P.~Manns, M.~Severitt, and S.~Strikwerda}, \emph{Input regularization for integer optimal control in {BV} with applications to control of poroelastic and poroviscoelastic systems}, J. Nonsmooth Anal. Optim., 5 (2024), 10529, \url{https://doi.org/10.46298/jnsao-2024-10529}.

\bibitem{BMW2023}
\textsc{L.~Bociu, B.~Muha, and J.~T. Webster}, \emph{Mathematical effects of linear visco-elasticity in quasi-static {B}iot models}, J. Math. Anal. Appl., 527 (2023), 127462, \url{https://doi.org/10.1016/j.jmaa.2023.127462}.

\bibitem{BS2022}
\textsc{L.~Bociu and S.~Strikwerda}, \emph{Optimal control in poroelasticity}, Appl. Anal., 101 (2022), pp.~1774--1796, \url{https://doi.org/10.1080/00036811.2021.2008372}.

\bibitem{BSvisco2022}
\textsc{L.~Bociu and S.~Strikwerda}, \emph{Poro-visco-elasticity in biomechanics: optimal control}, in Research in Mathematics of Materials Science, Assoc. Women Math. Ser.~31, Springer, Cham, 2022, pp.~103--132, \url{https://doi.org/10.1007/978-3-031-04496-0_5}.

\bibitem{ColliGilardiSprekels2018}
\textsc{P.~Colli, G.~Gilardi, and J.~Sprekels}, \emph{Optimal velocity control of a viscous {C}ahn--{H}illiard system with convection and dynamic boundary conditions}, SIAM J. Control Optim., 56 (2018), pp.~1665--1691, \url{https://doi.org/10.1137/17M1146786}.

\bibitem{Grisvard1985}
\textsc{P.~Grisvard}, \emph{Elliptic Problems in Nonsmooth Domains}, Monogr. Stud. Math.~24, Pitman, Boston, 1985.

\bibitem{GMV2026}
\textsc{T.~Gudi, D.~Mora, and N.~Verma}, \emph{Virtual element approximations for distributed or {N}eumann boundary optimal control problems governed by poroelasticity equation}, preprint, \url{https://arxiv.org/abs/2608.22107}, 2026.

\bibitem{Hinze2005}
\textsc{M.~Hinze}, \emph{A variational discretization concept in control constrained optimization: the linear-quadratic case}, Comput. Optim. Appl., 30 (2005), pp.~45--61, \url{https://doi.org/10.1007/s10589-005-4559-5}.

\bibitem{KLS2026}
\textsc{A.~Khan, J.~J. Lee, and H.~Singh}, \emph{Distributed optimal control problems governed by poroelasticity equations}, preprint, \url{https://arxiv.org/abs/2605.30839}, 2026.

\bibitem{KZ2026}
\textsc{C.~Kreuzer and P.~Zanotti}, \emph{Inf-sup theory for the quasi-static {B}iot's equations in poroelasticity}, SIAM J. Math. Anal., 58 (2026), pp.~1257--1287, \url{https://doi.org/10.1137/24M1674911}.

\bibitem{LeeLocking2026}
\textsc{J.~J. Lee}, \emph{Locking-free numerical methods for poroelasticity optimal control problems with minimal regularities}, preprint, \url{https://arxiv.org/abs/2611.22585}, 2026.

\bibitem{LeeThermo2026}
\textsc{J.~J. Lee}, \emph{Numerical approximation of distributed optimal control problems governed by linear thermo-poroelasticity}, preprint, \url{https://arxiv.org/abs/2609.24448}, 2026.

\bibitem{Lions1971}
\textsc{J.-L. Lions}, \emph{Optimal Control of Systems Governed by Partial Differential Equations}, Springer, Berlin, 1971.

\bibitem{LionsMagenes1972}
\textsc{J.-L. Lions and E.~Magenes}, \emph{Non-Homogeneous Boundary Value Problems and Applications}, Vol.~I, Grundlehren Math. Wiss.~181, Springer, Berlin, 1972, \url{https://doi.org/10.1007/978-3-642-65161-8}.

\bibitem{Showalter1974}
\textsc{R.~E. Showalter}, \emph{Degenerate evolution equations and applications}, Indiana Univ. Math. J., 23 (1974), pp.~655--677.

\bibitem{Showalter2000}
\textsc{R.~E. Showalter}, \emph{Diffusion in poro-elastic media}, J. Math. Anal. Appl., 251 (2000), pp.~310--340, \url{https://doi.org/10.1006/jmaa.2000.7048}.

\bibitem{Strikwerda2023}
\textsc{S.~Strikwerda}, \emph{Optimal Control in Fluid Flows through Deformable Porous Media}, Ph.D. thesis, North Carolina State University, 2023, \url{https://www.lib.ncsu.edu/resolver/1840.20/40780}.

\bibitem{Troeltzsch2010}
\textsc{F.~Tr\"oltzsch}, \emph{Optimal Control of Partial Differential Equations: Theory, Methods and Applications}, Grad. Stud. Math.~112, AMS, Providence, RI, 2010, \url{https://doi.org/10.1090/gsm/112}.

\bibitem{Verri2018}
\textsc{M.~Verri, G.~Guidoboni, L.~Bociu, and R.~Sacco}, \emph{The role of structural viscoelasticity in deformable porous media with incompressible constituents: applications in biomechanics}, Math. Biosci. Eng., 15 (2018), pp.~933--959, \url{https://doi.org/10.3934/mbe.2018042}.

\bibitem{Zenisek1984}
\textsc{A.~\v{Z}en\'{\i}\v{s}ek}, \emph{The existence and uniqueness theorem in {B}iot's consolidation theory}, Apl. Mat., 29 (1984), pp.~194--211, \url{https://doi.org/10.21136/AM.1984.104085}.
\end{thebibliography}
\end{document}